\documentclass[a4paper, 10pt]{article}

\usepackage[T1]{fontenc}
\usepackage[latin1]{inputenc}
\usepackage[english]{babel}
\usepackage{multicol}
\usepackage{fancyhdr}
\usepackage{amsmath,amssymb,amsthm}
\usepackage[all]{xy}
\usepackage{pstricks,pst-node}
\usepackage{geometry}
\renewcommand\c[1]{\mathcal{#1}}

\title{On Galois coverings and tilting modules}
\author{Patrick Le Meur
\footnote{\textit{adress:} Département de Mathématiques, Ecole normale
supérieure de Cachan, 61 avenue du président Wilson, 94235 Cachan, France}
\footnote{\textit{e-mail:} plemeur@dptmaths.ens-cachan.fr}}
\date{\today}

\newtheorem{prop}{Proposition}[section]
\newtheorem{Prop}[]{Proposition}
\newtheorem{thm}[prop]{Theorem}
\newtheorem{Thm}[]{Theorem}
\newtheorem{Cor}[]{Corollary}
\newtheorem{lem}[prop]{Lemma}
\newtheorem{cor}[prop]{Corollary}
\newtheorem{rem}[prop]{Remark}
\newtheorem{definition}[prop]{Definition}
\newtheorem{ex}[prop]{Example}

\def\c#1{\mathcal{#1}}
\def\e{\varepsilon}
\def\square{\blacksquare}

\makeatletter
\renewcommand\section{\@startsection {section}{1}{0mm}%
                                   {-3.5ex \@plus -1ex \@minus -.2ex}%
                                   {2.3ex \@plus.2ex}%
                                   {\normalfont\large\bfseries}}
\makeatother

\begin{document}
\maketitle
\abstract{
Let $A$ be a basic connected finite dimensional algebra over an
algebraically closed field $k$. Let $T$ be a basic tilting
$A$-module with arbitrary finite projective dimension. For a fixed group $G$ we compare
the set of isoclasses of connected Galois coverings of $A$ with group
$G$ and the set of isoclasses of connected Galois coverings of
$End_A(T)$ with group $G$. Using the Hasse diagram $\overrightarrow{\c K}_A$ (see
\cite{happel_unger} and \cite{riedtmann_schofield}) of basic
tilting $A$-modules, we give sufficient conditions on $T$ under which
there is a bijection between these two sets (these conditions are
always verified when $A$ is of finite representation type). Then we apply these
results to study when the simple connectedness of $A$ implies the one
of $End_A(T)$ (see \cite{assem_marcos_delapena}). Finally, using an
argument due to W. Crawley-Boevey, we prove that the type of any simply
connected tilted algebra is a tree and that its first Hochschild
cohomology group vanishes
}
\section*{Introduction}\label{intro}

Let $k$ be an algebraically closed field and let $A$ be a finite dimensional
$k$-algebra. In order to study the category $mod(A)$ of finite
dimensional (left) $A$-modules we may assume that $A$ is basic and
connected. In the study of $mod(A)$, tilting theory has proved to be a
powerful tool. Indeed, if $T$ is a basic tilting $A$-module and if we
set $B=End_A(T)$, then $A$ and $B$ have many common properties:
Brenner-Butler Theorem establishes an equivalence between certain
subcategories of $mod(A)$ and $mod(B)$ (see
\cite{brenner_butler}, \cite{happel_ringel} and \cite{miya}), $A$ and
$B$ have equivalent derived categories (see \cite{happel}) and
(in particular) they have isomorphic Grothendieck groups and
isomorphic Hochschild cohomologies.
In this text we will study the
following problem relating $A$ and $B$:
\begin{center}
\hfill\textit{is it possible to compare the Galois coverings of $A$ and those of
$B$?}\hfill$(P_1)$
\end{center}
As an example, if $A=kQ$ with $Q$ a finite quiver without oriented
cycle and if $T$ is an APR-tilting module associated to a sink $x$ of
$Q$ (see \cite{apr}) then
$B=kQ'$ where $Q'$ is obtained from $Q$ by reversing all the arrows
endings at $x$. In particular $Q$ and $Q'$ have the same underlying
graph and therefore $A$ has a connected Galois covering with group $G$
if and only if the same holds for $B$.

Recall that in order to consider Galois coverings of $A$ we always
consider $A$ as a $k$-category. When $\c C\to A$ is a Galois covering,
it is possible to describe part of 
$mod(A)$ in terms of $\c C$-modules (see for example
\cite{bongartz_gabriel} and \cite{gabriel}). This description is
useful because $mod(\c C)$ is easier to study than $mod(A)$, especially
when $\c C$ is simply connected (this last situation may occur when
$A$ is of finite representation type, see \cite{gabriel}).
Notice that
simple connectedness and tilting theory have already been studied
together through the following conjecture formulated in
\cite{assem_marcos_delapena}:
\begin{center}
\hfill$A$ is simply connected $\Longrightarrow$ $B$ is simply
connected\hfill $(P_2)$
\end{center}
More precisely, the above implication is true if: $A$ is of finite
representation type and $T$ is of projective dimension at most one
(see \cite{assem_skowronski2}), or if: $A=kQ$ (with $Q$ a quiver) and $B$ is
tame (see \cite{assem_marcos_delapena}, see also
\cite{assem_coelho_trepode} for a generalisation to the case of
quasi-tilted algebras). 
The two problems $(P_1)$ and $(P_2)$ are related
because $A$ is simply connected if and only if there is no proper
Galois covering $\c C\to A$ with $\c C$ connected and locally
bounded (see \cite{lemeur}). 

In order to study the question $(P_1)$ we will exhibit sufficient
conditions for $T$ to be of the first kind w.r.t. a fixed Galois
covering $\c C\xrightarrow{F} A$. Indeed, if $T$ is of the first kind
w.r.t. $F$, then it is possible to construct a Galois covering of
$B$. Under additional hypotheses on $T$, the equivalence class of
this Galois covering  is uniquely determined by the equivalence class
of $F$. Here we say that two Galois
coverings $F\colon\c C\to A$ and $F'\colon \c C'\to A$ are equivalent
if and only if there exists a commutative square of $k$-categories and
$k$-linear functors:
\begin{equation}
\xymatrix{
\c C\ar@{->}[d]_F\ar@{->}[r]^{\sim}&\c C'\ar@{->}[d]^{F'}\\
A\ar@{->}[r]^{\sim}&A
}\notag
\end{equation}
where horizontal arrows are isomorphisms and where the bottom
horizontal arrow restricts to the identity map on the set of objects
of $A$. For simplicity, let us say that $A$ and $B$ have
the same connected Galois coverings with group $G$ if there exists a
bijection between the sets $Gal_A(G)$ and $Gal_B(G)$ where $Gal_A(G)$
(resp. $Gal_B(G)$) stands for the set of equivalence classes of Galois
coverings $\c C\to A$ (resp. $\c C\to B$) with group $G$ and with $\c
C$ connected and locally bounded. With this definition, we
prove the following theorem which is the main result of this text and
which partially answers $(P_1)$:
\begin{Thm}
\label{thm0.1}
Let $T$ be a basic tilting $A$-module, let $B=End_A(T)$ and let $G$ be group.
\begin{enumerate}
\item If $T'\in \overrightarrow{\c K}_A$ lies in the connected component of
$\overrightarrow{\c K}_A$ containing $T$, then $End_A(T)$ and
$End_A(T')$ have the same connected Galois coverings with group $G$.
\item If $T$ lies in the connected component of $\overrightarrow{\c K}_A$
containing $A$ or $DA$, then $A$ and
$B$ have the same connected Galois coverings with group $G$.
\end{enumerate}
In particular, if $\overrightarrow{\c K}_A$ is connected (which
happens when $A$ is of finite representation type) then $A$ and $B$
have the same connected Galois coverings with group $G$, for any group $G$.
\end{Thm}
Here $\overrightarrow{\c K}_A$ is the Hasse diagram associated with the poset $\c
T_A$ of basic tilting $A$-modules (see \cite{happel_unger} and
\cite{riedtmann_schofield}). Recall (see \cite{gabriel}) that when $A$
is of finite representation type, $A$ admits a connected Galois covering with
group $G$ if and only if $G$ is a factor group of the fundamental
group $\pi_1(A)$ of the Auslander-Reiten quiver of $A$ with its mesh
relations. Theorem~\ref{thm0.1} allows us to get the
following corollary when $A$ and $B$ are of finite representation
type. We thank Ibrahim~Assem for having pointed out this corollary.
\begin{Cor}
\label{cor0.3}
Let $T$ be a basic tilting $A$-module and let $B=End_A(T)$. If both $A$ and $B$ are of finite
representation type, then $A$ and $B$ have isomorphic fundamental
groups.
\end{Cor}

Theorem~\ref{thm0.1} also allows us to prove the
following corollary related to $(P_2)$.
\begin{Cor}
\label{cor0.2}
(see \cite{assem_skowronski2} and \cite{assem_coelho_trepode})
Let $T$ be a basic tilting $A$-module and let $B=End_A(T)$.
\begin{enumerate}
\item If $T'\in\overrightarrow{\c K}_A$ lies in the connected component
of $\overrightarrow{\c K}_A$ containing $T$, then:
$End_A(T)$ is simply connected if and only if $End_A(T')$ is simply
  connected.
\item If $T$ lies in the connected component of
$\overrightarrow{\c K}_A$ containing $A$ or $DA$ then:
$A$ is simply connected if and only if $B$ is simply
  connected.
\end{enumerate}
In particular, f $\overrightarrow{\c K}_A$ is connected (e.g. $A$
is of finite representation type, see \cite{happel_unger}), then:
$A$ is simply connected if and only if $B$ is simply connected.
\end{Cor}

Theorem~\ref{thm0.1} uses the Hasse diagram
$\overrightarrow{\c K}_A$ in order to prove (in particular) that a
basic tilting $A$-module lying in a specific connected component of
$\overrightarrow{\c K}_A$ is of the first kind w.r.t. a given
connected Galois covering of $A$. On the other hand, an argument due to
W.~Crawley-Boevey (see
\cite{geiss_leclerc_schroer}) proves that any rigid $A$-module $M$
(i.e. such that $Ext^1_A(M,M)=0$) is of the first kind w.r.t. any
Galois covering of $A$ with group $\mathbb{Z}$ if $car(k)=0$
(resp. $\mathbb{Z}/p\mathbb{Z}$ if $car(k)=p$, $p$ prime). Using this
argument we are able to prove the following proposition. It gives a partial
answer to $(P_2)$ and to A.~Skwro\'nski's question in \cite[Problem
1]{skowronski2} wether it is true that a triangular algebra is simply
connected if and only if its first Hochschild cohomology group vanishes.
\begin{Prop}
\label{prop0.4}(see \cite{assem_coelho_trepode}, \cite{assem_marcos_delapena})
Assume that $A=kQ$ with $Q$ a quiver without oriented cycle. Let $T$
be a basic tilting $A$-module and let $B=End_A(T)$ (hence, $B$ is
tilted of type $Q$). Then the following implication holds:
\begin{center}
 $B$ is simply connected $\Longrightarrow$  $Q$ is
a tree and $HH^1(B) =0$
\end{center}
\end{Prop}
Recall that the implication of Proposition~\ref{prop0.4} has been proved
in \cite{assem_marcos_delapena} for $B$ tame tilted and in
\cite{assem_coelho_trepode} for $B$ tame quasi-tilted.

The text is organised as follows. In Section~\ref{section1} we will give the
definition of all the notions mentioned above and which will be used
for the proof of Theorem~\ref{thm0.1}. In Section~\ref{section2} we
will detail the construction
and give some properties of the Galois covering $F'$ of $B$ starting from a
Galois covering $F\colon \c C\to A$ of $A$ and a basic tilting
$A$-module $T$. In this
study, we will introduce the following hypotheses on the $A$-module $T$: $1$) $T$ is
of the first kind w.r.t. $F$ (this ensures that $F'$
 exists), $2$) the $\c C$-module $F.T$ obtained from $T$ by
restricting the scalars is basic (this ensures that $F'$ is connected
if $F$ is connected) $3$) $\psi.N\simeq
N$ for any direct summand $N$ of $T$ and for any automorphism
$\psi\colon A\xrightarrow{\sim} A$ which restricts to the identity map
on objects (this ensures that the equivalence class of $F'$ does
depend only on the equivalence class of $F$).
These three hypotheses lack of simplicity, therefore,
Section~\ref{section3} is devoted to find simple sufficient conditions
for the basic tilting $A$-module $T$ to verify these.
In particular, we will prove that the condition ``\textit{$T$ lies in
the connected component of 
$\overrightarrow{\c K}_A$ containing $A$}'' fits our 
requirements.
Since our main objective is to establish a correspondence between the
equivalence classes of the connected Galois coverings of $A$ and those
of $B$, we will
need to find conditions for $T$ to lie in both connected components
of $\overrightarrow{\c K}_A$ and $\overrightarrow{\c K}_B$ containing
$A$ and $B$ respectively (recall that $T$ is also a basic tilting $B$-module).
This will be done in Section~\ref{section4} where we compare the Hasse
diagrams $\overrightarrow{\c K}_A$ and $\overrightarrow{\c K}_B$. In
particular, we will prove that there is an oriented path in $\overrightarrow{\c
K}_A$ starting at $A$ and ending at $T$ if and only if there is an
oriented path in $\overrightarrow{\c K}_B$ starting at $B$ and ending
at $T$. This equivalence will be used in Section~\ref{section5} in
order to prove Theorem~\ref{thm0.1}, Corollary~\ref{cor0.3} and
Corollary~\ref{cor0.2}. Finally, in Section~\ref{sections}, we prove
Proposition~\ref{prop0.4}.

I would like to acknowledge Eduardo~N.~Marcos for his stimulating
remarks concerning the implication $(P_2)$ during the CIMPA school
\textit{Homological methods and representations of non-commutative
  algebras} in Mar del Plata, Argentina (February 2006).

\section{Basic definitions and preparatory lemmata}\label{section1}

\textbf{Reminder on $k$-categories} (see \cite{bongartz_gabriel} for
more details). A \textbf{$k$-category} is small category
$\c C$ such that for any $x,y\in Ob(\c C)$ the set $_y\c C_x$ of
morphisms from $x$ to $y$ is a $k$-vector space and such that the
composition of morphisms in $\c C$ is $k$-bilinear. A $k$-category $\c
C$ is called connected if and only if there is no non trivial partition $Ob(\c
C)=E\sqcup F$ such that $_y\c C_x=\ _x\c C_y=0$ for any $x\in E, y\in
F$.

All functors
between $k$-categories are supposed to be $k$-linear. If $F\colon \c
E\to \c B$ and $F'\colon \c E'\to \c B$ are functors between
$k$-categories, then $F$ and $F'$ are called equivalent if there exists a
commutative diagram:
\begin{equation}
\xymatrix{
\c E\ar@{->}[r]^{\sim} \ar@{->}[d]_F&\c E'\ar@{->}[d]^{F'}\\
\c B\ar@{->}[r]^{\sim}&\c B     
}\notag
\end{equation}
where horizontal arrows are isomorphisms and where the bottom
horizontal arrow restricts to the identity map on $Ob(\c B)$.
A \textbf{locally bounded} $k$-category is a $k$-category $\c C$
verifying the following conditions:
\begin{enumerate}
\item[.] distinct objects in $\c C$ are not isomorphic,
\item[.] for any $x\in Ob(\c C)$, the $k$-vector spaces
  $\bigoplus_{y\in Ob(\c C)}\ _y\c C_x$ and
$\bigoplus_{y\in Ob(\c C)}\ _x\c C_y$ are finite dimensional,
\item[.] for any $x\in Ob(\c C)$, the $k$-algebra $_x\c C_x$ is local.
\end{enumerate}
\noindent{}For example, let $A$ be a basic finite dimensional $k$-algebra
(basic means that $A$ is the direct sum of pairwise non-isomorphic
indecomposable projective $A$-modules) and let $\{e_1,\ldots,e_n\}$ be
a complete set of pairwise orthogonal primitive idempotents. Then $A$
can be viewed as a locally bounded $k$-category as follows:
$e_1,\ldots,e_n$ are the objects of $A$,
the space of morphisms from $e_i$ to $e_j$ is equal to
  $e_jAe_i$ for any $i,j$
and the composition of morphisms is induced by the product in
  $A$.
Notice that different choices for the primitive idempotents
$e_1,\ldots,e_n$ give rise to isomorphic $k$-categories. In this text
we shall always consider such an algebra $A$ as a locally bounded
$k$-category.\\

\textbf{Modules over $k$-categories}. If $\c C$ is a
$k$-category, a (left) $\c C$-module is a $k$-linear
functor $M\colon\c C\to MOD(k)$ where $MOD(k)$ is the category of
$k$-vector spaces. A morphism of $\c C$-modules $M\to N$ is a $k$-linear
natural transformation of functors. The category of $\c C$-modules is denoted
by $MOD(\c C)$.

A $\c C$-module $M$ is called locally finite
dimensional (resp. finite dimensional) if and only if $M(x)$ is finite dimensional
for any $x\in Ob(\c C)$ (resp. $\bigoplus_{x\in Ob(\c
  C)}M(x)$ is finite dimensional). The category of locally finite
dimensional (resp. finite dimensional) $\c C$-modules is denoted
by $Mod(\c C)$ (resp. $mod(\c C)$). Notice that if $\c C=A$ as above,
then $Mod(\c C)=mod(\c C)$.

We shall write $IND(\c C)$ (resp. $Ind(\c
C)$, resp. $ind(\c C)$) for the
full subcategory of $MOD(\c C)$ (resp. of $Mod(\c C)$, resp. of
$mod(\c C)$) of indecomposable $\c C$-modules. 
Finally, if
$M=N_1\bigoplus\ldots\bigoplus N_t$
with $N_i\in ind(\c C)$ for any $i$, then $M$ is
called \textbf{basic} if and only if $N_1,\ldots,N_t$ are pairwise
non isomorphic.\\

\textbf{Tilting modules}. Let $A$ be a basic finite dimensional 
$k$-algebra. A tilting $A$-module (see \cite{brenner_butler},
\cite{happel_ringel} and \cite{miya}) is a module $T\in mod(A)$
verifying the following conditions:
\begin{enumerate}
\item[$(T1)$] $T$ has finite projective dimension (i.e. $pd_A(T)<\infty$),
\item[$(T2)$] $Ext^i_A(T,T)=0$ for any $i>0$ (i.e. $T$ is selforthogonal),
\item[$(T3)$] there is an exact sequence in $mod(A)$: $0\to A\to
T_1\to\ldots\to T_r\to 0$ with $T_1,\ldots,T_r\in add(T)$
(this last property means that $T_1,\ldots,T_r$ are direct sums of direct summands of $T$).
\end{enumerate}
A module which satisfies conditions $(T1)$ and $(T2)$ above is
called an exceptional module.
Assume that $T$ is a tilting $A$-module. Then, $T$ is also a tilting $End_A(T)$-module
for the following action: $f.t=f(t)$ for $f\in End_A(T)$ and $t\in T$.
Assume moreover that $T$ is basic as an $A$-module and fix a decomposition
$T=T_1\oplus\ldots\oplus T_n$ with $T_1,\ldots,T_n\in ind(A)$. This
defines a decomposition of the unit of $End_A(T)$ into a
sum of primitive pairwise orthogonal idempotents so that $B:=End_A(T)$ is
a locally bounded $k$-category as follows: the set of objects is
$\{T_1,\ldots,T_n\}$ and for any $i,j$ the space of morphisms
$_{T_j}B_{T_i}$ is equal to $Hom_A(T_i,T_j)$. For any $x\in Ob(A)$,
$T(x)$ is an indecomposable $B$-module:
\begin{equation}
\begin{array}{rcl}
B & \rightarrow & MOD(k)\\
T_i\in Ob(B)&\mapsto &T_i(x)\\
u\in\ _{T_j}B_{T_i}&\mapsto&
T_i(x)\xrightarrow{u_x}T_j(x)
\end{array}\notag
\end{equation}
and $T=\bigoplus_{x\in Ob(A)}T(x)$. Finally, the following functor is an
isomorphism of $k$-categories:
\begin{equation}
\begin{array}{rrcl}
\rho_A\colon & A & \longrightarrow & End_B(T)\\
& x\in Ob(A) & \longmapsto & T(x)\in Ob(End_B(T))\\
& u\in\ _yA_x&\longmapsto & T(x)\xrightarrow{T(u)}T(y)
\end{array}\notag
\end{equation}
 For more details on the above properties and for a more general study
 of $End_A(T)$, we refer the reader to \cite{brenner_butler},
\cite{happel}, \cite{happel_ringel} and \cite{miya}.\\
Let $\c T_A$ be the set of basic tilting $A$-modules up to
isomorphism. Then $\c T_A$ is endowed with a partial order introduced
in \cite{riedtmann_schofield} and defined as follows. If $T\in\c T_A$,
the right perpendicular category $T^{\perp}$ of $T$ is defined by
(see \cite{auslander_reiten}):
\begin{equation}
T^{\perp}=\{X\in mod(A)\ |\ (\forall i\geqslant 1)\ \
Ext^i_A(T,X)=0\}\notag
\end{equation}
If $T'\in\c T_A$ is another basic tilting module, we write $T\leqslant
T'$ provided that $T^{\perp}\subseteq T^{'\perp}$. In particular, we
have $T\leqslant A$ for any $T\in\c T_A$. In \cite{happel_unger},
D.~Happel and L.~Unger  have proved that the Hasse diagram
$\overrightarrow{\c K}_A$ of $\c T_A$ is as follows. The vertices in
$\overrightarrow{\c K}_A$ are the elements in $\c T_A$ and there is an arrow
$T\to T'$ in $\overrightarrow{\c K}_A$ if and only if: $T=X\bigoplus \overline{T}$
with $X\in ind(A)$, $T'=Y\bigoplus \overline{T}$ with $Y\in ind(A)$ and there
exists a non split exact sequence $0\to X\xrightarrow{u} M\xrightarrow{v} Y\to 0$ in $mod(A)$
with $M\in add(\overline{T})$. In such a situation, $u$
(resp. $v$) is the left (resp. right) $add(\overline{T})$-approximation
of $X$ (resp. $Y$). For more details on
$\overrightarrow{\c K}_A$, we refer the reader to \cite{happel_unger} and
\cite{happel_unger2}.\\

\textbf{Galois coverings of $k$-categories}. Let $G$ be a group. A
\textbf{free $G$-category} is a $k$-category $\c E$ endowed with a
morphism of groups $G\to Aut(\c E)$ such that the induced action of
$G$ on $Ob(\c E)$ is free. In this case, there exists a (unique) quotient $\c
E\to \c E/G$ of $\c E$ by $G$ in the category of $k$-categories. With
this property, a Galois covering of $\c B$ with group $G$ is by
definition a functor $F\colon\c E\to \c B$ endowed with a group
morphism $G\to Aut(F)=\{g\in Aut(\c E)\ |\ F\circ g=F\}$ and verifying the
following facts:
\begin{enumerate}
\item[.] the group morphism $G\to Aut(F)\hookrightarrow Aut(\c E)$
  endows $\c E$ with a structure of free $G$-category,
\item[.] the functor $\c E/G\xrightarrow{\overline{F}}\c B$ induced by
  $F$ is an isomorphism.
\end{enumerate}
This definition implies that the group morphism $G\to
Aut(F)$ is one-to-one (actually one can show that this is an
isomorphism when $\c E$ is connected). Moreover for any $x\in Ob(\c
B)$ the set $F^{-1}(x)$ is non empty and called the fiber of $F$ at
$x$. It verifies  $F^{-1}(F(x))=G.x$ for any $x\in Ob(\c E)$.

We recall that Galois coverings are particular cases of \textbf{covering
functors} (see \cite{bongartz_gabriel}). A covering functor is a
$k$-linear functor $F\colon \c E\to\c B$ such that 
for any $x,y\in\c E_0$, the following mappings induced by $F$ are
bijective:
\begin{equation}
\bigoplus\limits_{y'\in F^{-1}(F(y))}\ _{y'}\c E_x\to \ _{F(y)}\c
B_{F(x)}\ \text{and}\ \bigoplus\limits_{x'\in F^{-1}(F(x))}\ _y\c E_{x'}\to \ _{F(y)}\c
B_{F(x)}\notag
\end{equation}
Remark that a covering functor is not supposed to restrict to a
surjective mapping on objects. However, a covering functor is an
isomorphism of $k$-categories if and only if it restricts to a
bijective mapping on objects. Using basic linear algebra arguments it
is easy to prove the following useful lemma:
\begin{lem}
\label{lem1.4}
Let $p,q$ be $k$-linear functors such that the composition $q\circ p$
is defined. Then $p,q$ and $q\circ p$ are covering functors as soon as
two of them are so.
\end{lem} 

If $F\colon \c E\to \c B$ is a Galois
covering with group $G$ and with $\c B$ connected then $\c E$ need not
be connected. In such a case, if
$\c E=\coprod\limits_{i\in I}\c E_i$
where the $\c E_i$'s are the connected components of $\c E$, then for
each $i$, the following functor:
\begin{equation}
F_i\colon \c E_i\hookrightarrow \c E\to\c B\notag
\end{equation}
is a Galois covering with group:
\begin{equation}
G_i:=\{g\in G\ |\ g(Ob(\c E_i))\cap Ob(\c E_i)\neq\emptyset\}=\{g\in
G\ |\ g(Ob(\c E_i))=Ob(\c E_i)\}\notag 
\end{equation}
Moreover, if $i,j\in I$ then the groups $G_i$ and $G_j$ are conjugated
in $G$ and there exists a commutative diagram:
\begin{equation}
\xymatrix{
\c E_i\ar@{->}[rr]^{\sim} \ar@{->}[rd]_{F_i} & & \c E_j
\ar@{->}[ld]^{F_j}\\
& \c B & 
}\notag
\end{equation}
where the horizontal arrow is an isomorphism. This implies 
that $G$ acts transitively on the set $\{\c E_i\ |\ i\in
I\}$ of the connected components of $\c E$. Notice that all these facts
may be false if $\c B$ is not connected.

Two Galois coverings of $\c B$ are called
equivalent if and only if they are isomorphic as functors between
$k$-categories (see above, this implies that the groups of the Galois
coverings are isomorphic). The equivalence class of a Galois covering
$F$ will be denoted by $[F]$. Finally, we shall say for short that \textbf{a
Galois covering $\c E\to \c B$ is connected} if and only if $\c E$ is
connected and locally bounded (this implies that $\c B$ is connected
and locally bounded, see \cite[1.2]{gabriel}).\\

\textbf{Simply connected locally bounded $k$-categories}. Let $\c B$
be a locally bounded $k$-category. Then $\c B$ is called
simply connected if and only if there is no proper connected Galois
covering of $\c B$ (proper means with non trivial group). This
definition is equivalent to the original one (see \cite{martinezvilla_delapena} for the
triangular case and \cite[Prop. 4.1]{lemeur} for the non-triangular
case) which was introduced in \cite{assem_skowronski}: $\c B$ is
simply connected if and only if $\pi_1(Q_{\c B},I)=1$ for any
admissible presentation $kQ_{\c B}/I\simeq \c B$ of $\c B$ (see
\cite{martinezvilla_delapena} for the definition of $\pi_1(Q_{\c
  B},I)$).\\

\textbf{Basic notions on covering techniques} (see \cite{bongartz_gabriel} and
\cite{riedtmann}). Let $F\colon \c E\to\c
B$ be a Galois covering with group $G$. The $G$-action on $\c E$ gives
rise to an action of $G$ on $MOD(\c
E)$: if $M\in MOD(\c E)$ and $g\in G$, then $^gM:=F\circ g^{-1}\in
MOD(\c E)$. Moreover, $F$ defines two
additive functors $F_{\lambda}\colon MOD(\c E)\to MOD(\c B)$ (the
push-down functor) and $F.\colon MOD(\c B)\to MOD(\c E)$ (the pull-up
functor) with the following properties
 (for more details we refer the reader
to \cite{bongartz_gabriel}):
\begin{enumerate}
\item[.] $F.M=M\circ F$ for any $M\in MOD(\c B)$,
\item[.] if $M\in MOD(\c E)$, then $\left(F_{\lambda}M\right)(x)=\bigoplus_{x'\in
    F^{-1}(x)}M(x')$ for any $x\in Ob(\c B)$. If $u\in\ _y\c E_x$,
  then the restriction of $\left(F_{\lambda}M\right)(F(u))$ to $M(g.x)$
(for $g.x\in F^{-1}(F(x))=G.x$) is equal to $M(g.u)\colon M(g.x)\to M(g.y)$,
\item[.] $F_{\lambda}$ and $F.$ are exact and send projective
  modules to projective modules,
\item[.] $F_{\lambda}\c E\simeq\bigoplus_{g\in G}\c B$ and $F.\c
B\simeq\c E$, where $\c E$ (resp. $\c B$) is the $\c E$-module
$x\mapsto \bigoplus_{y\in Ob(\c E)}\,_y\c E_x$ (resp. the $\c
B$-module $x\mapsto \bigoplus_{y\in Ob(\c B)}\,_y\c B_x$),
\item[.] $F.F_{\lambda}= \bigoplus_{g\in G}\ ^gId_{MOD(\c E)}$
\item[.] if $X\in MOD(\c B)$ verifies $X\simeq F_{\lambda}Y$ for
some $Y\in MOD(\c E)$, then $F_{\lambda}F.X\simeq \bigoplus_{g\in G} X$,
\item[.] $F_{\lambda}(mod(\c E))\subseteq mod(\c B)$,
  $F_{\lambda}(Mod(\c E))\subseteq Mod(\c B)$, $F.(Mod(\c
  B))\subseteq Mod(\c E)$,
\item[.] $D\circ F.=F.\circ D$ and $D\circ F_{\lambda|mod(\c
    E)}\simeq F_{\lambda}\circ D_{|mod(\c E)}$ where $D=Hom_k(?,k)$ is
  the usual duality,
\item[.] $F_{\lambda}$ is left adjoint to $F.$,
\item[.] $D\circ F_{\lambda}\circ D$ is right adjoint to $F.$ (in
  particular, there is a functorial isomorphism $Hom_{\c
    E}(F.M,N)\simeq Hom_{\c B}(M,F_{\lambda}N)$ for any $M\in
  MOD(\c B)$ and any $N\in mod(\c E)$).
\item[.] for any $M,N\in MOD(\c E)$, the following mappings induced by
  $F_{\lambda}$ are bijective:
 \begin{equation}
 \bigoplus\limits_{g\in G}Hom_{\c E}(\,^gM,N)\to Hom_{\c
   B}(F_{\lambda}M,F_{\lambda}N)\ \text{and}\ \bigoplus\limits_{g\in G}Hom_{\c E}(M,\,^gN)\to Hom_{\c
   B}(F_{\lambda}M,F_{\lambda}N)\notag
\end{equation}
\end{enumerate}
These properties give the following result which will be used many
times in this text:
\begin{lem}
\label{lem1.1}
If $M\in MOD(\c E)$ (resp. $M\in MOD(\c B)$) has finite projective
dimension, then so does $F_{\lambda}M$ (resp. $F.(M)$).

Let $M\in MOD(\c E)$, $N\in MOD(\c B)$ and $j\geqslant 1$. There is an
isomorphism of vector spaces:
\begin{equation}
Ext_{\c B}^j(F_{\lambda}M,N)\simeq Ext_{\c E}^j(M,F.N)\notag
\end{equation}
Moreover, if $M\in mod(\c E)$ then there is an isomorphism of vector spaces:
\begin{equation}
Ext_{\c E}^j(F.N,M)\simeq Ext_{\c B}^j(N,F_{\lambda}N)\notag
\end{equation}
\end{lem}
\noindent{\textbf{Proof:}} The first assertion is due to the fact that
$F.$ and $F_{\lambda}$ are exact and send projective modules to
projective modules. For
the same reasons, $F.$ and $F_{\lambda}$ induce $F.\colon\c
D(MOD(\c B))\to \c D(MOD(\c E))$ and $F_{\lambda}\colon \c D(MOD(\c
E))\to \c D(MOD(\c B))$ respectively and the adjunctions
$(F_{\lambda},F.)$ and $(F.,F_{\lambda})$ at the level of module
categories give rise to adjunctions at the level of derived
categories. Since $Ext_{\c E}^j(X,Y)=Hom_{\c D(MOD(\c E))}(Y,X[j])$ we
get the announced isomorphisms.\hfill$\square$\\

Remark that an isomorphism of $k$-categories is a particular case of
Galois covering. When $F$ is an isomorphism, $F.$ and $F_{\lambda}$
have additional properties as shows the following lemma whose proof
is a direct consequence of the definition of the push-down and pull-up
functors.
\begin{lem}
\label{lem1.2}
Assume that $F\colon \c E\to \c B$ is an isomorphism of $k$-categories. Then
$F.F_{\lambda}=Id_{MOD(\c E)}$ and $F_{\lambda}F.=Id_{MOD(\c
  B)}$.
\end{lem}

\textbf{Modules of the first kind}. Let $F\colon\c
E\to \c B$ be a Galois covering with group $G$. A  $\c
B$-module $M$ is called of the first kind w.r.t. $F$ if and only if
for any indecomposable direct summand $N$ of $M$ there exists
$\widehat{N}\in MOD(\c E)$ such that $N\simeq
F_{\lambda}\widehat{N}$. We will denote by $ind_1(\c B)$
(resp. $mod_1(\c B)$) the full subcategory of $ind(\c B)$ (resp. of
$mod(\c B)$) of modules of the first kind w.r.t. $F$. Notice the
following properties of $ind_1(\c B)$:
\begin{enumerate}
\item[.] if $M\in ind_1(\c B)$ and $N\in MOD(\c E)$ verify $M\simeq
F_{\lambda}N$, then $N\in ind(\c E)$,
\item[.] if $M\in ind_1(\c B)$ and $N,N'\in MOD(\c E)$ verify $M\simeq
F_{\lambda}N\simeq F_{\lambda}N'$, then there exists $g\in G$ such
that $N'\simeq\ ^gN$.
\end{enumerate}

If $\c B$ is connected and if $\c E=\coprod_{i\in I}\c E_i$, where the $\c
E_i$'s are the connected components of $\c E$, then an
indecomposable $\c B$-module $M$
is of the first kind w.r.t. $F$ if and only if it is
of the first kind w.r.t. $F_i\colon \c E_i\hookrightarrow \c E\to\c
B$ for any $i\in I$. More precisely, we have the following well know
lemma where we keep the established notations.
\begin{lem}
\label{lem1.3}
 Let $M\in ind(\c B)$. If $\widehat{M}\in ind(\c E)$ is such that
  $F_{\lambda}\widehat{M}\simeq M$, then there is a unique $i\in I$ such that
  $\widehat{M}\in ind(\c E_i)$. In such a case, we have $M\simeq
  (F_i)_{\lambda}\widehat{M}$. Moreover, if $j\in I$ then there
  exists $g\in G$ such that $g(\c E_i)=\c E_j$, and for any such $g$ we
  have: $^g\widehat{M}\in ind(\c E_j)$ and
  $(F_j)_{\lambda}\,^g\widehat{M}\simeq M$.
\end{lem}
\vskip 20pt
Throughout this text $A$ will denote a basic and
connected finite dimensional $k$-algebra and $n$ will denote the rank of its Grothendieck
group $K_0(A)$.

\section{Galois coverings associated with modules of the first kind}\label{section2}

 Throughout this section we will use the following data:
\begin{enumerate}
\item[-] $F\colon\c C\to A$  a Galois covering with
group $G$
\item[-] $T=T_1\bigoplus\ldots\bigoplus T_n\in
mod(A)$ (with $T_i\in ind(A)$)  a basic tilting $A$-module of the first kind
w.r.t. $F$,
\item[-] $\lambda_i\colon
F_{\lambda}(\widehat{T}_i)\to T_i$  an isomorphism with
$\widehat{T}_i\in ind(\c C)$, for every $i\in\{1,\ldots,n\}$.
\end{enumerate}
Let $B=End_A(T)$. With these data, we wish to:
\begin{enumerate}
\item construct a Galois covering $F_{\widehat{T}_i,\lambda_i}$ with
  group $G$ of $B$,
\item study the dependence of the equivalence class of $F_{\widehat{T}_i,\lambda_i}$
  on the data $\widehat{T}_i$, $\lambda_i$ and on the choice
  of $F$ in its equivalence class $[F]$,
\item repeat the construction made at the first step starting from $T$
  (viewed as a basic tilting $B$-module) and the Galois covering
  $F_{\widehat{T}_i,\lambda_i}$. This will give a Galois covering of
  $End_B(T)$ which will be compared with $F$ using
  the isomorphism $\rho_A\colon A\xrightarrow{\sim} End_B(T)$.
\end{enumerate}

\subsection{Construction of the Galois covering $F_{\widehat{T}_i,\lambda_i}$}
Let $\c End_{\c C}(\bigoplus_{g,i}\ ^g\widehat{T}_i)$ be the following
$k$-category:
\begin{enumerate}
\item[.] the set of objects is $\{\ ^g\widehat{T}_i\ |\ g\in G,\
  i\in\{1,\ldots,n\}\}$ ($^g\widehat{T}_i$ and $^{g'}\widehat{T}_j$
  are considered as different objects if $(i,g)\neq (j,g')$),
\item[.] the space of morphisms from $\,^g\widehat{T}_i$ to
  $\,^h\widehat{T}_j$ is equal to $Hom_{\c
    C}(\,^g\widehat{T}_i,\,^h\widehat{T}_j)$,
\item[.] the composition is induced by the composition of morphisms in
  $mod(\c C)$.
\end{enumerate}
\begin{rem}
\label{rem2.5}
\begin{enumerate}
\item The $\c C$-modules $\bigoplus_{g,i}\ ^g\widehat{T}_i$ and $F.T$ are
isomorphic.
\item If $G$ is a finite group, then $\c C$ is a finite dimensional
$k$-algebra. In particular, $\c End_{\c C}(\bigoplus_{g,i}\
^g\widehat{T}_i)$ and $End_{\c C}(F.T)$ are isomorphic $k$-algebras.
\item The $G$-action on $mod(\c C)$ naturally endows $\c End_{\c
C}(\bigoplus_{g,i}\,^g\widehat{T}_i)$ with a structure of free
$G$-category.
\item $\c End_{\c C}(\bigoplus_{i,g}\,^g\widehat{T}_i)$ is locally
bounded if and only if $G_{\widehat{T}_i}=1$ for any
$i$. This is equivalent to say that $F.T$ is a basic $\c C$-module.
\end{enumerate}
\end{rem}
The isomorphisms $\lambda_1,\ldots,\lambda_n$
define the following functor:
\begin{equation}
\begin{array}{cccc}
F_{\widehat{T}_i,\lambda_i}\colon&\c End_{\c C}(\bigoplus_{g,i}\ ^g\widehat{T}_i) & \longrightarrow & B\\
&^g\widehat{T}_i&\longmapsto & T_i\\
&^g\widehat{T}_i\xrightarrow{u}\,^h\widehat{T}_j&\longmapsto &
T_i\xrightarrow{\lambda_j\ \ F_{\lambda}u\ \ \lambda_i^{-1}}T_j
\end{array}\notag
\end{equation}

\begin{lem}
\label{lem2.9}
The functor $F_{\widehat{T}_i,\lambda_i}\colon \c
End_{\c C}(\bigoplus_{i,g}\,^g\widehat{T}_i)\to B$ is a Galois covering
with group $G$.
\end{lem}
\noindent{\textbf{Proof:}} For simplicity, we shall write $\c C'$ for $\c
End_{\c C}(\bigoplus_{i,g}\,^g\widehat{T}_i)$ and $F'\colon \c C'\to
B$ for $F_{\widehat{T}_i,\lambda_i}$. Recall (see
Remark~\ref{rem2.5}) that $G$ acts freely on $\c C'$. Moreover, we have
$F'\circ g=F'$ by construction of $F'$. So, $F'$ defines a commutative
diagram of $k$-categories and $k$-linear functors:
\begin{equation}
\xymatrix{
\c C'\ar@{->}[d] \ar@{->}[rd]^{F'}&\\
\c C'/G \ar@{->}[r]_{\overline{F'}} & B}\tag{$\star$}
\end{equation}
Where $\c C'\to\c C'/G$ is the quotient functor.
From the properties verified by
$F_{\lambda}$ (see Section~\ref{section1}) we infer that $F'$ is a
covering functor. Since $\c C'\to \c C'/G$ is
a also covering functor we deduce that so is $\overline{F'}$ (see
Lemma~\ref{lem1.4}). Finally, $\overline{F'}$ restricts to a bijective
mapping $Ob(\c C')/G=\{\,^g\widehat{T}_i\ |\ g\in G,\
i\in\{1,\ldots,n\}\}/G\to Ob(B)=\{T_1,\ldots,T_n\}$
so $\overline{F'}$ is an isomorphism. Thus, $F'$ is a Galois
covering with group $G$.\hfill$\square$\\

Since $F_{\widehat{T}_i,\lambda_i}$ is a Galois
covering, it is natural to ask
whether $\c End_{\c C}(\bigoplus_{i,g}\,^g\widehat{T}_i)$ is connected or not. The
following lemma partially answers this question.
\begin{lem}
\label{lem2.8}
If $\c C$ is not connected, then $\c End_{\c
    C}(\bigoplus_{i,g}\,^g\widehat{T}_i)$ is not connected.
\end{lem}
\noindent{\textbf{Proof}:} 
For simplicity let us write $\c C'$ for $\c End_{\c
  C}(\bigoplus_{i,g}\,^g\widehat{T}_i)$. Assume that $\c C$ is not
connected and let $\c
C=\coprod_{x\in I}\c C_x$ where the $\c C_x$'s are the connected
components of $\c C$. For $i\in\{1,\ldots,n\}$, we have
$\widehat{T}_i\in ind(\c C)$, so there exists a unique $x_i\in I$ such
that $\widehat{T}_i\in ind(\c C_{x_i})$. Let us set:
\begin{equation}
G_{x_1}=\{g\in G\ |\ g(\c C_{x_1})=\c C_{x_1}\}\notag
\end{equation}
Let $i\in\{1,\ldots,n\}$, since $G$ acts transitively on $\{\c C_{x}\
 |\ x\in I\}$,
 there exists $g_i\in G$ such that $g_i(\c
C_{x_1})=\c C_{x_i}$ (in particular $g_1\in G_{x_1}$).
Therefore:
\begin{equation}
(\forall i\in\{1,\ldots,n\})\ \ ^{g_i^{-1}}\widehat{T}_i\in mod(\c
C_{x_1})\notag
\end{equation}
Let us set $O$ to be the following set of objects of $\c C'$:
\begin{equation}
O:=\{\,^g\widehat{T}_i\ |\ i\in\{1,\ldots,n\}\ \text{and}\ gg_i\in
G_{x_1}\}\subseteq Ob(\c C')\notag
\end{equation}
Remark that $O$ satisfies the following:
\begin{enumerate}
\item[.] $O\neq\emptyset$ because $\widehat{T}_1\in O$.
\item[.] Since $\c C$
is not connected and since $G$ acts
transitively on $\{\c C_x\ |\ x\in I\}$ we have $G_{x_1}\subsetneq G$.
Let $g\in G\backslash G_{x_1}$, then $gg_1\not\in G_{x_1}$
 and $^g\widehat{T}_1\not\in O$. Hence
$O\subsetneq Ob(\c C)$.
\item[.] For any $^g\widehat{T}_i\in Ob(\c C')$, we have
  $^g\widehat{T}_i\in O$ if and only if $^g\widehat{T}_i\in ind(\c
  C_{x_1})$. As a consequence, there is no non zero morphism in $\c C'$
  between an object in $O$ and an object in $Ob(\c C')\backslash O$.
\end{enumerate}
As a consequence, $\c C'$ is not connected.\hfill$\square$\\

\subsection{Independence of the equivalence class of
$F_{\widehat{T}_i,\lambda_i}$ on the data $F, \widehat{T}_i,\lambda_i$}
In the two following lemmas, we examine the dependence of the
equivalence class $[F_{\widehat{T}_i,\lambda_i}]$ of 
$F_{\widehat{T}_i,\lambda_i}$ on the choice of
$\widehat{T}_1,\ldots,\widehat{T}_n,\lambda_1,\ldots,\lambda_n$ and on the
choice of $F$ in its equivalence class $[F]$.
\begin{lem}
\label{lem2.2}
For each $i\in\{1,\ldots,n\}$, let $\mu_i\colon
F_{\lambda}\overline{T}_i\to T_i$ be an isomorphism with
$\overline{T}_i\in ind(\c C)$. Then $F_{\widehat{T}_i,\lambda_i}$ and
$F_{\overline{T}_i,\mu_i}$ are equivalent.
\end{lem}
\noindent{\textbf{Proof:}} We need to exhibit a commutative square:
\begin{equation}
\xymatrix{
\c End_{\c C}(\bigoplus_{i,g}\,^g\overline{T}_i)\ar@{->}[r]^{\varphi}
\ar@{->}[d]_{F_{\overline{T}_i,\mu_i}} & \c End_{\c
  C}(\bigoplus_{i,g}\,^g\widehat{T}_i)
\ar@{->}[d]^{F_{\widehat{T}_i,\lambda_i}}\\
B\ar@{->}[r]^{\psi}&B
}\tag{$\star$}
\end{equation}
where $\varphi,\psi$ are isomorphisms and where $\psi(x)=x$ for any
$x\in Ob(B)=\{T_1,\ldots,T_n\}$. Let $i\in\{1,\ldots,n\}$.
We have $F_{\lambda}\overline{T}_i\simeq
T_i\simeq F_{\lambda}\widehat{T}_i$, so there exists an isomorphism
$\theta_i\colon\overline{T}_i\xrightarrow{\sim}\,^{g_i}\widehat{T}_i$
with $g_i\in G$. Let us define $\varphi$ by:
\begin{equation}
\begin{array}{rccc}
\varphi\colon&\c End_{\c C}(\bigoplus_{i,g}\,^g\overline{T}_i) & \longrightarrow & \c End_{\c
  C}(\bigoplus_{i,g}\,^g\widehat{T}_i)\\
&^g\overline{T}_i&\longmapsto&\,^{gg_i}\widehat{T}_i\\
&^g\overline{T}_i\xrightarrow{u}\,^h\overline{T}_j&\longmapsto&\,^{gg_i}\widehat{T}_i\xrightarrow{\,^h\theta_j\
  \,u\ \ ^g\theta_i^{-1}} \,^{hg_j}\widehat{T}_j
\end{array}
\notag
\end{equation}
Then $\varphi$ is an isomorphism of $k$-categories. Notice that
$\theta_i$ defines an isomorphism $F_{\lambda}\theta_i\colon
F_{\lambda}\overline{T}_i\to F_{\lambda}\widehat{T}_i$. So we can
define $\psi$ by:
\begin{equation}
\begin{array}{rccc}
\psi\colon&B&\longrightarrow&B\\
&T_i&\longmapsto & T_i\\
&T_i\xrightarrow{u}T_j&\longmapsto & \psi(u)
\end{array}\notag
\end{equation}
where $\psi(u)$ is the composition:
\begin{equation}
T_i \xrightarrow{\lambda_i^{-1}} F_{\lambda}\widehat{T}_i
\xrightarrow{F_{\lambda}\theta_i^{-1}} F_{\lambda}\overline{T}_i
\xrightarrow{\mu_i} T_i \xrightarrow{u} T_j \xrightarrow{\mu_j^{-1}}
F_{\lambda}\overline{T}_j \xrightarrow{F_{\lambda}\theta_j}
F_{\lambda}\widehat{T}_j\xrightarrow{\lambda_j} T_j\notag
\end{equation}
So $\psi$ is an isomorphism of $k$-categories which restricts to
the identity map on $Ob(B)$. Moreover $\varphi$ and $\psi$
make $(\star)$ commutative.\hfill$\square$\\

In the following lemma, we show that, under additional hypotheses on
$T$, the equivalence class of $F_{\widehat{T}_i,\lambda_i}$ does not
depend on the choice of $F$ in $[F]$.
\begin{lem}
\label{lem2.3}
Assume that $F'\colon\c C'\to A$ is a Galois covering (with group $G$)
equivalent to $F$ and assume that
$T$ verifies the following condition:
\begin{center}
($H_{A,T}$):''$\psi.T_i\simeq T_i$ for any $i$ and for any isomorphism $\psi\colon
A\xrightarrow{\sim} A$ which restricts to the identity map on $Ob(A)$.''
\end{center}
Then $T$ is of the first
 kind w.r.t. $F'$.
 For each 
$i\in\{1,\ldots,n\}$ let $\mu_i\colon F'_{\lambda}\overline{T}_i\to
T_i$ be an isomorphism with $\overline{T}_i\in ind(\c C')$. Then
$F'_{\overline{T}_i,\mu_i}$ and $F_{\widehat{T}_i,\lambda_i}$ are equivalent.
\end{lem}
\noindent{\textbf{Proof:}} Let us fix an isomorphism between $F$ and
 $F'$:
\begin{equation}
\xymatrix{
\c C'\ar@{->}[r]^{\varphi}\ar@{->}[d]_{F'} & \c C\ar@{->}[d]^F\\
A\ar@{->}[r]^{\psi}&A
}\notag
\end{equation}
Let  us set $\nu\colon Aut(\c C')\to Aut(\c
C)$ to be the isomorphism
of groups (recall that $Aut(\c C')=G$ and $Aut(\c C)=G$):
\begin{equation}
\begin{array}{crcl}
\nu\colon & Aut(\c C') & \to & Aut(\c C)\\
& g & \mapsto & \varphi\circ g\circ\varphi^{-1}
\end{array}\notag
\end{equation}
Recall that any $g\in Aut(\c C)=G$ (resp. $g\in Aut(\c C')=G$) defines
an automorphism $g$ of $MOD(\c C)$ (resp. of
$MOD(\c C')$). Therefore we have an equality of functors $MOD(\c
C')\to MOD(A)$:
\begin{equation}
(\forall g\in Aut(\c C'))\ \ \varphi_{\lambda}\circ g=\nu(g)\circ \varphi_{\lambda}\notag
\end{equation}

Let us fix an isomorphism $\theta_i\colon
\psi.T_i\to T_i$, for each $i$, and let us set
$\overline{T}_i=\varphi.\widehat{T}_i$. In particular:
$\varphi_{\lambda}\overline{T}_i=\widehat{T}_i$ (see
Lemma~\ref{lem1.2}). Since $\psi.\psi_{\lambda}=Id_{MOD(A)}$
(loc. cit.) and $\psi F'=F\varphi$, we infer that:
\begin{equation}
F'_{\lambda}\overline{T}_i = \psi.\psi_{\lambda}F'_{\lambda}\overline{T_i} =
\psi_.F_{\lambda}\varphi_{\lambda}\overline{T}_i = \psi.F_{\lambda}\widehat{T}_i\notag
\end{equation}
Therefore, we get for each $i$ an isomorphism $\mu_i\colon
F'_{\lambda}\overline{T}_i\to T_i$ equal to the composition:
\begin{equation}
\mu_i\colon F'_{\lambda}\overline{T}_i =
\psi_.F_{\lambda}\widehat{T}_i \xrightarrow{\psi.\lambda_i}
\psi_.T_i \xrightarrow{\theta_i} T_i\notag
\end{equation}
This proves that $T$ is of the first kind w.r.t. $F'$. According to
the preceding subsection, this defines the Galois covering with group
$G$:
\begin{equation}
F'_{\overline{T}_i,\mu_i}\colon \c End_{\c 
  C'}(\bigoplus_{g,i}\,^g\overline{T}_i) \to B\notag
\end{equation}
 Thanks to
Lemma~\ref{lem2.2} we only need to prove that $F'_{\overline{T}_i,\mu_i}$
and $F_{\widehat{T}_i,\lambda_i}$ are equivalent.

First, we have the following functor induced by $\varphi_{\lambda}$:
\begin{equation}
\begin{array}{rccc}
\overline{\varphi}\colon & \c End_{\c C'}(\bigoplus_{i,g}\,^g\overline{T}_i)
& \longrightarrow & \c End_{\c C}(\bigoplus_{i,g}\,^g\widehat{T}_i)\\
&^g\overline{T}_i & \longmapsto &
\,^{\nu(g)}\widehat{T}_i=\varphi_{\lambda}\,^g\overline{T}_i \\
 & ^g\overline{T}_i\xrightarrow{u}\,^h\overline{T}_j & \longmapsto &
 \,^{\nu(g)}\widehat{T}_i\xrightarrow{\varphi_{\lambda}u}\,^{\nu(h)}\widehat{T}_j
\end{array}\notag
\end{equation}
Since $\nu\colon G\to G$ is an isomorphism and because of the equalities
$\varphi_{\lambda}\varphi.=Id_{MOD(\c C)}$ and
$\varphi.\varphi_{\lambda}=Id_{MOD(\c C)}$ (see Lemma~\ref{lem1.2}),
the functor $\overline{\varphi}$ is an isomorphism.

Secondly, we have the following functor induced by $\psi_{\lambda}$:
\begin{equation}
\begin{array}{rccc}
\overline{\psi}\colon & B &\longrightarrow & B\\
& T_i & \longmapsto & T_i\\
& T_i\xrightarrow{u}T_j & \longmapsto &
T_i\xrightarrow{\psi_{\lambda}(\theta_{j}^{-1}u\theta_i)} T_j
\end{array}\notag
\end{equation}
Since $\psi_{\lambda}\psi.=\psi.\psi_{\lambda}=Id_{MOD(A)}$, the
functor $\overline{\psi}$ is a well defined isomorphism and restricts
to the identity map on
$Ob(B)$. Therefore, we have a diagram whose horizontal arrows
are isomorphisms and whose bottom horizontal arrow restricts to the
identity map on the set of objects:
\begin{equation}
\xymatrix{
\c End_{\c C'}(\bigoplus_{i,g}\,^g\overline{T}_i)
\ar@{->}[r]^{\overline{\varphi}} \ar@{->}[d]_{F'_{\overline{T}_i,\mu_i}} &
\c End_{\c C}(\bigoplus_{i,g}\,^g\widehat{T}_i)
\ar@{->}[d]^{F_{\widehat{T}_i,\lambda_i}}\\
B\ar@{->}[r]^{\overline{\psi}} & B
}\notag
\end{equation}
This diagram is commutative, indeed, for any
$^g\overline{T}_i\xrightarrow{u}\,^h\overline{T}_j$ we have:
\begin{equation}
\begin{array}{rclcl}
\overline{\psi}F'_{\overline{T}_i,\mu_i}(u) & = &
\overline{\psi}(\mu_jF'_{\lambda}(u)\mu_i^{-1})  = 
\psi_{\lambda}(\theta_j^{-1}\mu_jF'_{\lambda}(u)\mu_i^{-1}\theta_i)&&\\
&=&
\psi_{\lambda}(\theta_j^{-1}\theta_j\psi.(\lambda_j)F'_{\lambda}(u)\psi.(\lambda_i)^{-1}\theta_i^{-1}\theta_i)
& & \\
&=&
\lambda_j\ (\psi_{\lambda}F'_{\lambda})(u)\ \lambda_i^{-1} & &
\text{because $\psi_{\lambda}\psi.=Id_{MOD(\c C)}$}\\
&=&
\lambda_j\ (F_{\lambda}\varphi_{\lambda})(u)\ \lambda_i^{-1} & &
\text{because $F\varphi=\psi F'$}\\
&=&F_{\widehat{T}_i,\lambda_i}(\varphi_{\lambda}(u))=F_{\widehat{T}_i,\lambda_i}\overline{\varphi}(u)&&
\end{array}\notag
\end{equation}
This proves that $F'_{\overline{T}_i,\mu_i}$ and $F_{\widehat{T}_i,\lambda_i}$
are equivalent.\hfill$\square$\\

\label{rem2.6}
Later, we shall prove that if $T$ is a basic tilting $A$-module lying
in the connected component of $\overrightarrow{\c K}_A$ containing
$A$, then the hypothesis ($H_{A,T}$) in the preceding lemma is
automatically verified. As a consequence, for
these tilting $A$-modules, the property ``to be of the
first kind w.r.t. an equivalence class of Galois coverings of $A$''
does make sense. However, the hypothesis ($H_{A,T}$) is not
verified for any $A$-module $T$ as the following example shows.

\begin{ex}
\label{ex2.6}
Let $A$ be the path algebra
of the following quiver:
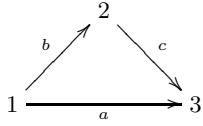
\begin{equation}
\xymatrix{
& 2\ar@{->}[rd]^c&\\
1\ar@{->}[ru]^b \ar@{->}[rr]_a&&3
}\notag
\end{equation}
Here $n=3$ and we have an isomorphism of $k$-categories:
$\psi\colon  A \xrightarrow{\sim}  A$ such that $\psi(x)=x$ for any
$x\in Ob(A)$, $\psi(a)=a+cb$, $\psi(b)=b$ and $\psi(c)=c$.
For any integer $i$, let $T_i$ be the following $A$-module:
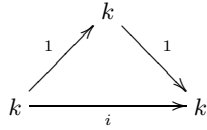
\begin{equation}
\xymatrix{
& k\ar@{->}[rd]^1&\\
k\ar@{->}[ru]^1 \ar@{->}[rr]_i&&k
}\notag
\end{equation}
Then:
\begin{enumerate}
\item[-] $T_i$ and $T_{i+1}$ are not isomorphic, for any $i$,
\item[-] if $car(k)\neq 2$, then $T_1,T_2,T_3$ are pairwise non
isomorphic,
\item[-] $\psi.T_i=T_{i+1}$ for any $i$.
\end{enumerate}
In particular, if $car(k)\neq 2$, then hypothesis ($H_{A,T}$) is not
satisfied for $T=T_1\oplus T_2\oplus T_3$. Remark that $T$ is
not tilting. Indeed, for any $i$, we have $Ext^1_A(T_i,T_i)\simeq k$
because $\tau_A(T_i)\simeq T_i$.
\end{ex}

Remark~\ref{rem2.5}, Lemma~\ref{lem2.2} and Lemma~\ref{lem2.3} justify
the following definition:
\begin{definition}
\label{def2.4}
Assume that the hypothesis ($H_{A,T}$) is satisfied (see Lemma~\ref{lem2.3}).
The equivalence class $[F]$ of $F$ and the basic tilting $A$-module
$T=T_1\bigoplus\ldots\bigoplus T_n$ of the first kind w.r.t. $[F]$ (with
$T_i\in ind(A)$) uniquely define an
equivalence class of Galois covering of $B$ with group $G$ and which admits
$F_{\widehat{T}_i,\lambda_i}$ as a representative. This
equivalence class will be denoted by $[F]_T\colon \c End_{\c C}(F.T)\to
B$.
\end{definition}

\subsection{Comparison of $[F]$ and $\left([F]_T\right)_T$}
For short, let us write $\c C'$
for $\c End_{\c
  C}(\bigoplus_{g,i}\,^g\widehat{T}_i)$ and $F'$ for
$F_{\widehat{T}_i,\lambda_i}$ .
In this subsection we shall not assume that the hypothesis ($H_{A,T}$)
of Lemma~\ref{lem2.3} is satisfied, except for the last proposition.
Starting from $F$ and from the isomorphisms $\lambda_i\colon
F_{\lambda}\widehat{T}_i\xrightarrow{\sim} T_i$, $i\in\{1,\ldots,n\}$,
we have constructed the Galois covering $F'$ of
$B$. One may try to perform the same construction
starting
from $F'$ in order to get a
Galois covering $F''$ of $End_B(T)\simeq A$ and eventually compare $F''$ with
$F$. In this purpose, we need to prove that $T$ is of
the first kind w.r.t. $F'$. Let us fix a lifting
$L\colon Ob(A)\to Ob(\c C)$ of the surjective mapping $F\colon Ob(\c
C)\to Ob(A)$. For $x\in Ob(A)$, let
$\widehat{T(x)}$ be the $\c C'$-module such that:
\begin{enumerate}
\item[-] $\widehat{T(x)}(\,^g\widehat{T}_i)=\widehat{T}_i(g^{-1}L(x))$
  for any $^g\widehat{T}_i\in Ob(\c C')$,
\item[-]
  $\widehat{T(x)}\left(\,^g\widehat{T}_i\xrightarrow{u}\,^h\widehat{T}_j\right)$
 is equal to
 $\widehat{T}_i(g^{-1}L(x))\xrightarrow{u_{L(x)}}\widehat{T}_j(h^{-1}L(x))$ 
for any $u\in\ _{^h\widehat{T}_j}\c C'_{^g\widehat{T}_i}$. 
\end{enumerate}
Therefore, for any $i\in\{1,\ldots,n\}$, we have:
\begin{equation}
\left(F'_{\lambda}\widehat{T(x)}\right)(T_i)=\bigoplus\limits_{g\in
  G}\widehat{T(x)}(\,^g\widehat{T}_i)=\bigoplus\limits_{g\in
  G}\widehat{T}_i(g^{-1}L(x))=\left(F_{\lambda}\widehat{T}_i\right)(x)\notag
\end{equation}
So we may set $(\mu_x)_{T_i}\colon
\left(F'_{\lambda}\widehat{T(x)}\right)(T_i)\to \left(T(x)\right)(T_i)$ to be equal
to $(F_{\lambda}\widehat{T_i})(x)\xrightarrow{(\lambda_i)_x}
T_i(x)$.
\begin{lem}
\label{lem2.1}
The linear isomorphisms $(\mu_x)_{T_i}$ ($i\in\{1,\ldots,n\}$) define
an isomorphism of $B$-modules:
\begin{equation}
\mu_x\colon F'_{\lambda}\widehat{T(x)}\xrightarrow{\sim} T(x)\notag
\end{equation}
\end{lem}
\noindent{\textbf{Proof:}} We only need to prove that $\mu_x$ is a
morphism of $B$-modules. Let $u\in\ _{^t\widehat{T}_j}\c
C'_{^s\widehat{T}_i}$ so that $F'(u)\in\ _{T_j}B_{T_i}$, and let us prove that the
following diagram commutes:
\begin{equation}
  \xymatrix{
\left(F'_{\lambda}\widehat{T(x)}\right)(T_i)
\ar@{->}[rrr]^{(\mu_x)_{T_i}=(\lambda_i)_x}
\ar@{->}[d]_{\left(F'_{\lambda}\widehat{T(x)}\right)(F'(u))} &&&
T_i(x) \ar@{->}[d]^{\left(T(x)\right)(F'(u))=F'(u)_x}\\
\left(F'_{\lambda}\widehat{T(x)}\right)(T_j)
\ar@{->}[rrr]^{(\mu_x)_{T_j}=(\lambda_j)_x} &&&
T_j(x)
}\tag{$\star$}
\end{equation}
Let $g\in G$ and let us compute the restriction of $F'(u)_x\circ (\lambda_i)_x$
to $\widehat{T(x)}(\,^g\widehat{T}_i)$. Recall that $F'(u)_x$ is
equal to the composition:
\begin{equation}
T_i(x)\xrightarrow{(\lambda_i^{-1})_x}\left(F_{\lambda}\,^s\widehat{T}_i\right)(x)\xrightarrow{\left(F_{\lambda}u\right)_x}
  \left(F_{\lambda}\,^t\widehat{T}_j\right)(x)
    \xrightarrow{(\lambda_j)_x}T_j(x)\notag
\end{equation}
Moreover, the restriction of
$\left(F_{\lambda}u\right)_x$ to
    $\widehat{T(x)}(\,^g\widehat{T}_i)=\widehat{T}_i(g^{-1}L((x))$ is (by construction of the push-down
    functor) equal to
    $\widehat{T}_i(g^{-1}L(x))\xrightarrow{u_{sg^{-1}L(x)}}\widehat{T}_j(t^{-1}sg^{-1}L(x))$.
Thus, the restriction of $F'(u)_x\circ (\lambda_i)_x$ to
$\widehat{T(x)}(\,^g\widehat{T}_i)$ is equal to the composition:
\begin{equation}
\widehat{T}_i(g^{-1}L(x))\xrightarrow{u_{sg^{-1}L(x)}}\widehat{T}_j(t^{-1}sg^{-1}L(x))\xrightarrow{(\lambda_j)_x}
T_j(x)\tag{$i$}
\end{equation}
On the other hand, the restriction of
$\left(F'_{\lambda}\widehat{T(x)}\right)(F'(u))$ to
$\widehat{T(x)}(\,^g\widehat{T}_i)=\widehat{T(x)}(\,^{gs^{-1}s}\widehat{T}_i)$
is (by construction of the 
push-down functor) equal to:
\begin{equation}
  \widehat{T(x)}(\,^g\widehat{T}_i)\xrightarrow{\widehat{T(x)}(\,^{gs^{-1}}u)}
  \widehat{T(x)}(\,^{gs^{-1}t}\widehat{T}_j)\tag{$ii$}
\end{equation}
and
$\widehat{T(x)}(\,^{gs^{-1}}u)=\left(\,^{gs^{-1}}u\right)_{L(x)}=u_{sg^{-1}L(x)}$.
This last equality, together with ($i$) and ($ii$), proves that the
diagram ($\star$) commutes.\hfill$\square$\\

Thanks to Lemma~\ref{lem2.1}, we have a Galois covering
$F'_{\widehat{T(x)},\mu_x}\colon \c End_{\c
  C'}(\bigoplus_{g,x}\,^g\widehat{T(x)})\to End_B(T)$ with group
$G$. For short, we shall write 
$\c C''$ for $\c End_{\c
  C'}(\bigoplus_{g,x}\,^g\widehat{T(x)})$ and $F''$ for
$F'_{\widehat{T(x)},\mu_x}$. The following lemma relates
$F''$ and $F$.
\begin{lem}
\label{lem2.10}
There exists an isomorphism of $k$-categories $\psi\colon\c
C\xrightarrow{\sim} \c C''$ such that the following diagram is
commutative:
\begin{equation}
\xymatrix{
\c C \ar@{->}[d]_F \ar@{->}[r]^{\psi} & \c C'' \ar@{->}[d]^{F''}\\
A \ar@{->}[r]^{\rho_A} & End_B(T)
}\notag
\end{equation}
In particular, $F$ and $\rho_A^{-1}F''$ are equivalent as Galois
coverings of $A$.
\end{lem}
\noindent{\textbf{Proof:}} Since $G$ acts freely on $Ob(\c C)$ and
since $L\colon Ob(A)\to Ob(\c C)$ lifts $F\colon Ob(\c C)\to Ob(A)$,
any $x\in Ob(\c C)$ is equal to $gL(x')$ with $g\in G,x'\in Ob(A)$
uniquely determined by $x$. Let $\psi\colon \c C\to \c C''$ be as
follows:
\begin{enumerate}
\item[-]
$\psi(gL(x))=\,^g\widehat{T(x)}$ for any $gL(x)\in Ob(\c C)$,
\item[-] for $u\in\ _{hL(y)}\c C_{gL(x)}$, we let
  $\psi(u)\colon\,^g\widehat{T(x)}\to \,^h\widehat{T(y)}$ be the
  morphism of $\c C'$-modules such that for any $^s\widehat{T}_i\in
  Ob(\c C')$, $\psi(u)_{\,^s\widehat{T}_i}$ is equal to:
\begin{equation}
\widehat{T}_i(s^{-1}u)\colon \widehat{T}_i(s^{-1}gL(x))\to
\widehat{T}_i(s^{-1}hL(y))\notag
\end{equation}
\end{enumerate}
Let us prove the following facts:
\begin{enumerate}
\item $\psi(u)$ is a morphism of $\c C'$-modules for any $u\in\
  _{hL(y)}\c C_{gL(x)}$,
\item $\psi$ is a functor,
\item $F''\circ\psi=\rho_A\circ F$,
\item $\psi$ is an isomorphism.
\end{enumerate}
$1$) Let $u\in\
  _{hL(y)}\c C_{gL(x)}$. We need to prove that for any $f\in\
  _{^t\widehat{T}_j}\c C'_{^s\widehat{T}_i}$, the following diagram
  commutes:
\begin{equation}
\xymatrix{
^g\widehat{T(x)}(\,^s\widehat{T}_i)
\ar@{->}[rrr]^{\psi(u)_{^s\widehat{T}_i}}
\ar@{->}[d]_{^g\widehat{T(x)}(f)} &&&
^h\widehat{T(y)}(\,^s\widehat{T}_i) \ar@{->}[d]^{^h\widehat{T(y)}(f)}
\\
^g\widehat{T(x)}(\,^t\widehat{T}_j)
\ar@{->}[rrr]^{\psi(u)_{^t\widehat{T}_j}}&& & 
^h\widehat{T(y)}(\,^t\widehat{T}_j)
}\notag
\end{equation}
By construction, this diagram is equal to:
\begin{equation}
\xymatrix{
\widehat{T}_i(s^{-1}gL(x)) \ar@{->}[rrr]^{\widehat{T}_i(s^{-1}u)}
\ar@{->}[d]_{f_{gL(x)}} &&&
\widehat{T}_i(s^{-1}hL(y)) \ar@{->}[d]^{f_{hL(y)}}\\
\widehat{T}_j(t^{-1}gL(x)) \ar@{->}[rrr]^{\widehat{T}_j(t^{-1}u)} &&&
\widehat{T}_j(t^{-1}hL(y))
}\notag
\end{equation}
and the latter is commutative because $f\colon\,^s\widehat{T}_i\to
\,^t\widehat{T}_j$ is a morphism of $\c C$-modules. This proves that
$\psi(u)$ is a morphism of $\c C'$-modules for any morphism $u$ in $\c C$. 

$2$) One easily checks that $\psi(1_{gL(x)})=Id_{^g\widehat{T(x)}}$ for
any $gL(x)\in Ob(\c C)$. Let $u,v$ be morphisms in $\c C$ such that
the composition $vu$ exists. Then, for any $^s\widehat{T}_j\in Ob(\c C')$:
\begin{equation}
\left(\psi(v)\circ
  \psi(u)\right)_{^s\widehat{T}_j}=\psi(v)_{^s\widehat{T}_j}\circ
\psi(u)_{^s\widehat{T}_j}=\widehat{T}_j(s^{-1}v)\circ
\widehat{T}_j(s^{-1}u)=\widehat{T}_j(s^{-1}(v\circ u))=\psi(v\circ
u)_{^s\widehat{T}_j}\notag
\end{equation}
So $\psi\colon \c C\to \c C''$ is a functor.

$3)$ Let $gL(x)\in Ob(\c C)$. Then:
\begin{equation}
F''\circ
\psi(gL(x))=F''(\,^g\widehat{T(x)})=T(x)=\rho_A(x)=\rho_A\circ
F(gL(x))\notag
\end{equation}
Let $u\in\ _{hL(y)}\c C_{gL(x)}$ and let us prove that $F''
\psi(u)=\rho_A F(u)$. Let $T_i\in Ob(B)$. Then:
\begin{align}
\left\{
\begin{array}{rcl}
\left(F''
\psi(u)\right)_{T_i} & = &T_i(x)\xrightarrow{(\lambda_i^{-1})_x}\left(F'_{\lambda}\,^g\widehat{T(x)}\right)(T_i)
  \xrightarrow{\left(F'_{\lambda}(\psi(u))\right)_{T_i}}
  \left(F'_{\lambda}\,^h\widehat{T(y)}\right)(T_i)
\xrightarrow{(\lambda_i)_y}
  T_i(y)\\
\left(\rho_A F(u)\right)_{T_i} &= &T_i(x)\xrightarrow{T_i(F(u))}T_i(y)
\end{array}
\right.\notag
\end{align}
Recall that
$\left(F'_{\lambda}\,^g\widehat{T(x)}\right)(T_i)=\bigoplus_{s\in
G}\,^g\widehat{T(x)}(\,^s\widehat{T}_i)$ and that
$^g\widehat{T(x)}(\,^s\widehat{T}_i)=\widehat{T}_i(s^{-1}gL(x))$, for
any $s\in G$. Let $s\in G$. Then, the restriction of
$\left(F'_{\lambda}(\psi(u))\right)_{T_i}$ to
$^g\widehat{T(x)}(\,^s\widehat{T}_i)$ is equal to:
\begin{equation}
\widehat{T}_i(s^{-1}gL(x))
\xrightarrow{\psi(u)_{^s\widehat{T}_i}=\widehat{T}_i(s^{-1}u)}\widehat{T}_i(s^{-1}hL(y))\notag
\end{equation}
Therefore, for any $s\in G$, 
 the restriction of
  $(\lambda_i^{-1})_y\circ\left(F''
\psi(u)\right)_{T_i}\circ (\lambda_i)_x$ (resp.
$(\lambda_i^{-1})_y\circ\left(\rho_A F(u)\right)_{T_i}\circ
(\lambda_i)_x$) to $^g\widehat{T(x)}(\,^s\widehat{T}_i)=\widehat{T}_i(s^{-1}gL(x))$ is equal to:
\begin{equation}
\begin{array}{rl}
&\widehat{T}_i(s^{-1}gL(x))\xrightarrow{\ \ \widehat{T}_i(s^{-1}u)\ \ } \widehat{T}_i(s^{-1}hL(y))
\\
\text{resp.}&
\widehat{T}_i(s^{-1}gL(x))\xrightarrow{\ \ (\lambda_i^{-1})_y\
  T_i(F(u))\ (\lambda_i)_x\ \ } \widehat{T}_i(s^{-1}hL(y))
\end{array}\notag
\end{equation}
Since $\lambda_i\colon F_{\lambda}\left(\,^s\widehat{T}_i\right)\to T_i$ is an isomorphism of
$A$-modules,  $(\lambda_i^{-1})_y\circ T_i(F(u))\circ
(\lambda_i)_x$ equals $\widehat{T}_i(s^{-1}u)$. We infer that
$(\lambda_i^{-1})_y\circ\left(F'' 
\psi(u)\right)_{T_i}\circ (\lambda_i)_x$ and
$(\lambda_i^{-1})_y\circ\left(\rho_A F(u)\right)_{T_i}\circ
(\lambda_i)_x$ coincide on
$\widehat{T}(s^{-1}gL(x))$, for any $s\in G$. As a consequence, $\left(F''
\psi(u)\right)_{T_i}=\left(\rho_A F(u)\right)_{T_i}$, for any $T_i\in
Ob(B)$. This proves that $F''\circ
\psi(u)=\rho_A\circ F(u)$ for any morphism $u$ in $\c C$. In other
words: $F''\circ \psi=\rho_A\circ F$.

$4$) Let us prove that $\psi$ is an isomorphism. Since $F''$ and
$\rho_A\circ F$ are covering functors, Lemma~\ref{lem1.4} implies that
so does $\psi$. Since $\psi$ restricts to a
bijective mapping on objects, we deduce that $\psi$ is an isomorphism.
\hfill$\square$\\

Thanks to Lemma~\ref{lem2.10} we can complete Lemma~\ref{lem2.8}
concerning the connectedness of $\c C'$. The following proposition
will be useful in the sequel, it is a direct consequence of
Lemma~\ref{lem2.8} and Lemma~\ref{lem2.10}.
\begin{prop}
\label{prop2.11}
 $\c C$ is connected if and only if $\c C'=\c End_{\c
  C}(\bigoplus_{g,i}\widehat{T}_i)$ is connected.
\end{prop}

We finish this subsection with the following proposition which compares
the equivalence class of $F$ and $([F]_T)_T$ when the latter is well
defined (see Definition~\ref{def2.4}). It is a direct consequence of
Lemma~\ref{lem2.10}. Notice that $\rho_A^{-1}\circ ([F]_T)_T$ is an
equivalence class of Galois coverings of $A$.
\begin{prop}
\label{prop2.7}
Assume that both conditions ($H_{A,T}$) and ($H_{B,T}$) are
satisfied. Then, the equivalence class $[F]$ of $F$ coincides with
$\rho_A^{-1}\circ ([F]_T)_T$.
\end{prop}

\section{Tilting modules of the first kind}\label{section3}

Let $F\colon\c C\to A$ be a Galois covering with group
$G$ and with $\c C$ locally bounded. The aim of this section is to
give ``simple'' sufficient conditions which guarantee the following facts:
\begin{enumerate}
\item[-] $T$ is of the first kind w.r.t. $F$,
\item[-] $F.T$ is a basic
$\c C$-module,
\item[-] the hypothesis ($H_{A,T}$) is satisfied (see
  Lemma~\ref{lem2.3}), i.e. $\psi.N\simeq N$ for any direct summand
  $N$ of $T$ and for any automorphism $\psi\colon A\xrightarrow{\sim}
  A$ which restricts to the identity map on objects.
\end{enumerate}
We begin with the following proposition.
\begin{prop}
\label{prop5.1}
\label{prop5.6}
Assume that $T$ and $T'$ lie in a same connected component of $\overrightarrow{\c K}_A$. Then:
\begin{equation}
T\in mod_1(A)\Leftrightarrow T'\in mod_1(A)\notag
\end{equation}
In particular, if $T'=A$ or $T'=DA$, then $T\in mod_1(A)$.
\end{prop}
\noindent{\textbf{Proof:}} Since $A,DA\in mod_1(A)$, we only need to
prove the equivalence of the proposition under the assumption: there
is an arrow $T\to T'$ in $\overrightarrow{\c K}_A$. Let us assume that $T\in mod_1(A)$. Since
$T\to T'$ is an arrow in $\overrightarrow{\c K}_A$, we have the following
data:
\begin{enumerate}
\item[.] $T=X\bigoplus \overline{T}$ with $X\in ind(A)$,
\item[.] $T'=Y\bigoplus \overline{T}$ with $Y\in ind(A)$,
\item[.] $\e\colon 0\to X\to M\to Y\to 0$ a non split exact sequence
  in $mod(A)$ with $M\in add(\overline{T})$.
\end{enumerate}
Thus, we only need to prove that $Y\in mod_1(A)$ in order
to get $T'\in mod_1(A)$. In this purpose, we will need the
following lemma.
\begin{lem}
\label{lem5.3}
Let $\e:0\to X\xrightarrow{u} M\to Y\to 0$
be an exact sequence in $mod(A)$ verifying the following hypotheses:
\begin{enumerate}
\item[.] $X,Y\in ind(A)$ and $X=F_{\lambda}\widehat{X}$ (with
$\widehat{X}\in ind(\c C)$),
\item[.] $M=M_1\bigoplus \ldots \bigoplus M_t$ where
$M_i=F_{\lambda}\widehat{M}_i\in ind(A)$ (with $\widehat{M}_i\in
  ind(\c C)$), for every $i$,
\item[.] $Ext_A^1(Y,M)=0$.
\end{enumerate}
Then $(\e)$ is isomorphic to an exact sequence in $mod(A)$:
\begin{equation}
0\to X \xrightarrow{\begin{bmatrix} F_{\lambda}u_1'\\ \vdots \\
    F_{\lambda}u_t'\end{bmatrix}} M_1\bigoplus\ldots\bigoplus
    M_T\to Y\to 0\notag
\end{equation}
where $u'_i\in Hom_{\c C}(\widehat{X},\,^{g_i}\widehat{M}_i)$ for some
$g_i\in G$, for every $i$.
\end{lem}
\noindent{\textbf{Proof of Lemma~\ref{lem5.3}}}: 
For short, we shall say that $u\in Hom_A(X,M_i)$ is homogeneous of
degree $g\in G$ if and only if $u=F_{\lambda}u'$ with $u'\in
Hom_{\c C}(\widehat{X},\,^g\widehat{M}_i)$. Recall from
Section~\ref{section1} that any $u\in Hom_A(X,M_i)$ is (uniquely) the
sum of $d$ homogeneous morphisms with pairwise different degrees (with
$d\geqslant 0$).
Let us write
    $u=\begin{bmatrix} u_1\\ \vdots \\
    u_t\end{bmatrix}$ with $u_i\colon X\to M_i$ for each $i$. We may
  assume that $u_1\colon X\to M_1$ is not
    homogeneous. Thus:
\begin{equation}
u_1=h_1+\ldots+h_d\notag
\end{equation}
where $d\geqslant 2$ and $h_1,\ldots,h_d\colon X\to M_1$ are non zero
homogeneous morphisms of pairwise different degree. In order to prove
the lemma, it
suffices to prove the following property which we denote by $(\c P)$:
\begin{center}
``
\textit{
$(\e)$ is isomorphic to an exact sequence of the form:
\begin{equation}
0\to X\xrightarrow{\begin{bmatrix} u_1'\\ u_2\\ \vdots \\
    u_t\end{bmatrix}} M_1\bigoplus\ldots\bigoplus M_t\to Y\to 0
    \tag{$\e'$}
\end{equation}
where $u_1'$ is the sum of at most $d-1$ non zero homogeneous
morphisms of pairwise different degree.
}
``
\end{center}
For simplicity we adopt the following notations:
\begin{enumerate}
\item[.] $\overline{M}=M_2\bigoplus\cdots\bigoplus M_t$ (so
  $M=M_1\bigoplus \overline{M}$),
\item[.] $\overline{u}=\begin{bmatrix} u_2\\ \vdots \\
    u_t\end{bmatrix}\colon X\to \overline{M}$ (so $u=\begin{bmatrix}
    u_1\\\overline{u} \end{bmatrix}\colon X\to M_1\bigoplus\overline{M}$),
\item[.] $\overline{h}=h_2+\ldots+h_d\colon X\to M_1$ (so
  $u_1=h_1+\overline{h}$).
\end{enumerate}
From $Hom_A(\e,M_1)$ we get the exact sequence:
\begin{equation}
Hom_A(M_1\bigoplus\overline{M},M_1)\xrightarrow{u^*} Hom_A(X,M_1) \to
Ext_A^1(Y,M_1)=0\notag
\end{equation}
So there exists $[\lambda,\mu]\colon M_1\bigoplus\overline{M}\to M_1$ such
that $h_1=[\lambda,\mu]u$.
Hence:
\begin{equation}
h_1=\lambda u_1+\mu\overline{u}=\lambda
h_1+\lambda\overline{h}+\mu\overline{u}\tag{$i$}
\end{equation}
Let us distinguish two cases whether $\lambda\in End_A(M_1)$ is
invertible or nilpotent (recall that $M_1\in ind(A)$):\\
$\bullet$ If $\lambda$ is invertible then:
\begin{equation}
\theta:=\begin{bmatrix}
\lambda & \mu\\
0 & Id_{\overline{M}}
\end{bmatrix}\colon M_1\bigoplus \overline{M}
\to M_1\bigoplus \overline{M}\notag
\end{equation}
is invertible. Using $(i)$ we deduce an isomorphism of exact
sequences:
\begin{equation}
\xymatrix{
0\ar@{->}[r] & X \ar@{->}[r]^{\begin{bmatrix} u_1\\
    \overline{u}\end{bmatrix}} \ar@{=}[d] & M_1\bigoplus \overline{M}
\ar@{->}[r] \ar@{->}[d]^{\theta} & Y \ar@{->}[r]
\ar@{->}[d]^{\sim}& 0 & (\e)\\
0 \ar@{->}[r] & X \ar@{->}[r]^{\begin{bmatrix} h_1\\
    \overline{u}\end{bmatrix}} & M_1\bigoplus\overline{M}\ar@{->}[r] &
Y\ar@{->}[r] & 0 & (\e')
}\notag
\end{equation}
Since $h_1\colon X\to M_1$ is homogeneous, $(\e')$ fits property $(\c
P)$. So $(\c P)$ is satisfied in this case.\\
$\bullet$ If $\lambda\in End_A(M_1)$ is nilpotent, let $p\geqslant 0$
be such that $\lambda^p=0$. Using $(i)$ we get the following
equalities:
\begin{equation}
\begin{array}{rcl}
h_1&=&\lambda^2h_1+(\lambda^2+\lambda)\overline{h}+(\lambda+Id_{M_1})\mu\overline{u}\\
\vdots&\vdots&\vdots\\
h_1&=&\lambda^th_1+(\lambda^t+\lambda^{t-1}+\ldots+\lambda)\overline{h}+(\lambda^{t-1}+\ldots+\lambda+Id_{M_1})\mu\overline{u}\\
\vdots&\vdots&\vdots\\
h_1&=&\lambda^ph_1+(\lambda^p+\lambda^{p-1}+\ldots+\lambda)\overline{h}+(\lambda^{p-1}+\ldots+\lambda+Id_{M_1})\mu\overline{u}
\end{array}\notag
\end{equation}
Since $\lambda^p=0$ and $u_1=h_1+\overline{h}$ we infer that:
\begin{equation}
u_1=\lambda'\overline{h}+\lambda'\mu\overline{u}\notag
\end{equation}
where $\lambda':=Id_{M_1}+\lambda+\ldots+\lambda^{p-1}\in End_A(M_1)$ is
invertible. So we have an isomorphism:
\begin{equation}
\theta:=\begin{bmatrix}
\lambda'&\lambda'\mu\\
0&Id_{\overline{M}}
\end{bmatrix}\colon M_1\bigoplus\overline{M}\to M_1\bigoplus\overline{M}\notag
\end{equation}
and consequently we have an isomorphism of exact sequences:
\begin{equation}
\xymatrix{
0\ar@{->}[r] & X \ar@{=}[d]^{}
 \ar@{->}[r]^{\begin{bmatrix}
\overline{h}\\
\overline{u}
\end{bmatrix}} &M_1\bigoplus\overline{M} \ar@{->}[d]^{\theta}
\ar@{->}[r]& Y\ar@{->}[r] \ar@{->}[d]^{\sim}& 0 & (\e')\\
0\ar@{->}[r] & X\ar@{->}[r]^{
\begin{bmatrix}
u_1\\
\overline{u}
\end{bmatrix}
} & M_1\bigoplus \overline{M} \ar@{->}[r] & Y \ar@{->}[r] & 0 & (\e)
}\notag
\end{equation}
where $\overline{h}=h_2+\ldots+h_p$ is the sum of $p-1$ non zero homogeneous
morphisms of pairwise different degrees.
So $(\c P)$ is satisfied in this case. This finishes the proof
of the lemma.\hfill$\square$\\

Now we can prove that $Y\in mod_1(A)$. Thanks to the preceding
lemma, and with the same notations, we know that $(\e)$ is isomorphic to an exact sequence in
$mod(A)$:
\begin{equation}
0\to X\xrightarrow{\begin{bmatrix}
F_{\lambda}u_1'\\ \vdots \\ F_{\lambda}u_t'
\end{bmatrix}
} M_1\bigoplus\ldots M_t\to Y\to 0\tag{$\e'$}
\end{equation}
where $u'_i\in Hom_{\c C}(\widehat{X},\,^{g_i}\widehat{M}_i)$ for some
$g_i\in G$, for every $i$.
Therefore (recall that $F_{\lambda}$ is exact):
\begin{equation}
Y\simeq Coker\begin{bmatrix}
F_{\lambda}(u_1')\\ \vdots\\ F_{\lambda}(u_t')
\end{bmatrix}
\simeq F_{\lambda}\left(Coker
\begin{bmatrix}
u_1'\\ \vdots\\ u_t'
\end{bmatrix}
\right)\notag
\end{equation}
This proves that $Y\in mod_1(A)$. Therefore $T'=Y\bigoplus
\overline{T}\in mod_1(A)$.

The proof of the implication $T'\in mod_1(A)\Rightarrow T\in mod_1(A)$ is
similar, excepted that instead of using
Lemma~\ref{lem5.3} we use a dual version:
\begin{lem}
\label{lem5.4}
Let $\e\colon 0\to X\to M\xrightarrow{v} Y\to 0$
be an exact sequence in $mod(A)$ verifying the following hypotheses:
\begin{enumerate}
\item[.] $X,Y\in ind(A)$ and $Y=F_{\lambda}\widehat{Y}$ (with
$\widehat{Y}\in ind(\c C)$),
\item[.] $M=M_1\bigoplus \ldots \bigoplus M_t$ where
$M_i=F_{\lambda}\widehat{M}_i\in ind(A)$ (with $\widehat{M}_i\in
ind(\c C)$), for every $i$,
\item[.] $Ext_A^1(M,X)=0$.
\end{enumerate}
Then $(\e)$ is isomorphic to an exact sequence in $mod(A)$:
\begin{equation}
0\to X\to  M_1\bigoplus\ldots\bigoplus
    M_T\xrightarrow{\begin{bmatrix} F_{\lambda}v_1'\\ \vdots \\
    F_{\lambda}v_t'\end{bmatrix}} Y\to 0\notag
\end{equation}
where $v'_i\in Hom_{\c C}(\,^{g_i}\widehat{M}_i,\widehat{Y})$ for some
$g_i\in G$, for every $i$.
\end{lem}
This finishes the proof of Proposition~\ref{prop5.1}.\hfill$\square$\\

\begin{rem}
\label{rem5.7}
Proposition~\ref{prop5.1} is similar to part of \cite[Thm
3.6]{gabriel} where P.~Gabriel proves the following: if $F\colon \c
C\to A$ is a Galois covering with group $G$, with $\c C$ locally
bounded and such that $G$ acts freely on $ind(\c C)$, then for any
connected component $C$ of the Auslander-Reiten quiver of $A$, all
indecomposable modules of $C$ lie in $ind_1(A)$ as soon as any one of
them does.
\end{rem}

\begin{rem}
\label{rem5.5}
The proof of Proposition~\ref{prop5.1} shows that for an arrow $T\to
T'$ in $\overrightarrow{\c K}_A$ such that $T,T'\in mod_1(A)$ there
exists an exact sequence in $mod(\c C)$:
\begin{equation}
0\to X\xrightarrow{\iota} M\xrightarrow{\pi} Y\to 0\notag
\end{equation}
with the following properties:
\begin{enumerate}
\item[.] $T=F_{\lambda}X\bigoplus \overline{T}$ and
$F_{\lambda}X\in ind(A)$,
\item[.] $T'=F_{\lambda}Y\bigoplus \overline{T}$ and
$F_{\lambda}Y\in ind(A)$,
\item[.] $F_{\lambda}M\in add(\overline{T})$.
\end{enumerate}
\end{rem}

Recall that $\overrightarrow{\c K}_A$ has a Brauer-Thrall type
property (see \cite[Cor. 2.2]{happel_unger}): 
$\overrightarrow{\c K}_A$ is finite and connected if it has a finite
connected component. In particular, $\overrightarrow{\c K}_A$ is
finite and connected if $A$ is of finite representation type. Using
Proposition~\ref{prop5.6}, we get the following corollary.
\begin{cor}
\label{cor5.8}
If $\overrightarrow{\c K}_A$ is finite (e.g. $A$ is of finite
representation type), then any
$T\in\overrightarrow{\c K}_A$ is of the first kind w.r.t. $F$.
\end{cor}

Now we turn to the second goal of this section: for $T\in mod_1(A)$ a
basic tilting $A$-module, give sufficient conditions for $F.T$ to be a
basic $\c C$-module.
\begin{prop}
\label{prop5.9}
Let $T,T'\in \overrightarrow{\c K}_A\cap mod_1(A)$ lie in a same connected
component of $\overrightarrow{\c K}_A$, then:
\begin{center}
$F.T$ is a basic $\c C$-module $\Leftrightarrow$ $F.T'$ is a basic $\c C$-module.
\end{center}
In particular, if $T'=A$ or $T'=DA$, then $T\in mod_1(A)$ and $F.T$ is
a basic $\c C$-module.
\end{prop}
\noindent{\textbf{Proof}}: The $k$-category $\c C$ is locally bounded
so $F.A\simeq \c C$ and $F.(DA)\simeq D\c C$ are basic $\c
C$-modules.Therefore, we only need to prove the equivalence of the
proposition.Without loss of generality, we may assume that
there is an arrow $T\to T'$ in $\overrightarrow{\c K}_A$. Let us
assume that $F.T$ is basic and let us prove that so is $F.T'$.  We will
use Remark~\ref{rem5.5} from which the adopt the notations, in
particular, the exact sequence $0\to X\xrightarrow{\iota}
M\xrightarrow{\pi} Y\to 0$ in $mod(\c C)$ will be denoted by $(\e)$. 
Because $F.T$ is basic and because of the properties verified by
$(\e)$, we only need to prove that $G_{Y}=1$. 
Let $\varphi\colon Y\to\,^gY$ be an isomorphism in $mod(\c C)$ (with
$g\in G$), and
let us prove that $g=1$. To do this we will exhibit an isomorphism
$\theta\colon X\to\,^gX$. Notice that:
\begin{equation}
(\forall h\in G)\ \ \left\{
\begin{array}{l}
^hX,\,^hM\in add(F.T)\\
^hY,\,^hM\in add(F.T')
\end{array}\right.
\end{equation}
Moreover, thanks to $T\in\overrightarrow{\c K}_A$ and to
$F_{\lambda}F.T=\bigoplus_{h\in G}T$ (which is true because $T$ is of
the first kind w.r.t. $F$, see Section~\ref{section1}), we have:
\begin{equation}
(\forall i\geqslant 1)\ \ Ext_{\c C}^i(F.T,F.T)\simeq
Ext_A^i(F_{\lambda}F.(T),T)\simeq \prod\limits_{h\in
G}Ext_A^i(T,T)=0\notag
\end{equation}
In particular:
\begin{equation}
Ext_{\c C}^1(\,^gM,X)=Ext_{\c C}^1(M,\,^gX)=0\tag{$i$}
\end{equation}
With $Hom_{\c C}(M,\,^g\e)$, this last equality gives the exact
sequence:
\begin{equation}
Hom_{\c C}(M,\,^gM)\xrightarrow{(\,^g\pi)_*}Hom_{\c C}(M,\,^gY)\to
Ext_{\c C}^1(M,\,^gX)=0\notag
\end{equation}
From this exact sequence, we deduce the existence of $\psi\in Hom_{\c
C}(M,\,^gM)$ such that the following diagram commutes:
\begin{equation}
\xymatrix{
M\ar@{->}[r]^{\pi} \ar@{->}[d]_{\psi} & Y\ar@{->}[d]_{\varphi}\\
^gM\ar@{->}[r]^{\,^g\pi}&^gY
}\notag
\end{equation}
This implies the existence of $\theta\in Hom_{\c C}(X,\,^gX)$ making
commutative the following diagram with exact rows:
\begin{equation}
\xymatrix{
0\ar@{->}[r] & X\ar@{->}[r]^{\iota} \ar@{->}[d]_{\theta} &
M\ar@{->}[r]^{\pi} \ar@{->}[d]_{\psi} & Y \ar@{->}[d]_{\varphi}
\ar@{->}[r] & 0\\
0\ar@{->}[r] & ^gX\ar@{->}[r]^{^g\iota} &
^gM\ar@{->}[r]^{^g\pi} & ^gY 
\ar@{->}[r] & 0
}\tag{$ii$}
\end{equation}
We claim that $\theta\colon X\to\,^gX$ is an isomorphism. The
arguments that have been used to get $(ii)$ may be adapted (just
replace the use of $Hom_{\c C}(M,\,^g\e)$ and of $\varphi\colon
Y\to\,^gY$ by $Hom_{\c C}(\,^gM,\e)$ and $\varphi^{-1}\colon \,^gY\to
Y$) to get the following commutative diagram with exact rows:
\begin{equation}
\xymatrix{
0\ar@{->}[r] & ^gX\ar@{->}[r]^{^g\iota} \ar@{->}[d]_{\theta'} &
^gM\ar@{->}[r]^{^g\pi} \ar@{->}[d]_{\psi'} & ^gY \ar@{->}[d]_{\varphi^{-1}}
\ar@{->}[r] & 0\\
0\ar@{->}[r] & X\ar@{->}[r]^{\iota} &
M\ar@{->}[r]^{\pi} & Y 
\ar@{->}[r] & 0
}\tag{$iii$}
\end{equation}
In order to show that $\theta\colon X\to\, ^gX$ is an isomorphism, let
us show that $\theta'\theta\in End_{\c C}(X)$ is an
isomorphism. Notice $(ii)$ and $(iii)$ give the following commutative
diagram:
\begin{equation}
\xymatrix{
0\ar@{->}[rr] && X\ar@{->}[rr]^{\iota} \ar@{->}[d]_{\theta'\theta-id_X} &&
M\ar@{->}[rr]^{\pi} \ar@{->}[d]_{\psi'\psi-Id_M} && Y \ar@{->}[d]_{0}
\ar@{->}[rr] && 0\\
0\ar@{->}[rr] && X\ar@{->}[rr]^{\iota} &&
M\ar@{->}[rr]^{\pi} && Y 
\ar@{->}[rr] && 0
}\tag{$iv$}
\end{equation}
In particular we have $\pi(\psi'\psi-Id_M)=0$, so there exists
$\lambda\in Hom_{\c C}(M,X)$ such that:
\begin{equation}
\psi'\psi-Id_M=\iota\lambda\notag
\end{equation}
Therefore:
\begin{equation}
\iota(\theta'\theta-Id_X)=\iota\lambda\iota\notag
\end{equation}
Since $\iota$ is one-to-one, we get $\theta'\theta-id_X=\lambda\iota$,
i.e.:
\begin{equation}
\theta'\theta=Id_X+\lambda\iota\notag
\end{equation}
If $\lambda\iota\in End_{\c C}(X)$ was an isomorphism, then
$\iota\colon X\to M$ would be a section. This would imply that $F_{\lambda}X$
is a direct summand of $F_{\lambda}M$. This last property is
impossible because: $T=F_{\lambda}X\bigoplus\overline{T}$,
$F_{\lambda}M\in add(\overline{T})$ and $T$ is basic. This
contradiction proves that $\lambda\iota\in End_{\c C}(X)$ is
nilpotent. Therefore $\theta'\theta=Id_X+\lambda\iota\in End_{\c
C}(X)$ is invertible. As a consequence, $\theta\colon X\to \,^gX$ is a
section. Since $X,\,^gX\in ind(\c C)$, we deduce that $\theta\colon
X\to \,^gX$ is an isomorphism. But we assumed that $F.T$ is basic,
so $g=1$. This finishes the proof of the implication:
\begin{center}
$F.T$ is basic $\Rightarrow$ $F.T'$ is basic.
\end{center}
After exchanging the roles of $T$ and $T'$ in the above
arguments, we also prove that:
\begin{center}
$F.T$ is basic $\Leftarrow$ $F.T'$ is basic.
\end{center}
under the assumption that $T\to T'$ is an arrow in $\overrightarrow{\c K}_A$.
This achieves the proof of the proposition.\hfill$\square$\\

Proposition~\ref{prop5.9} has the following corollary which will be useful
to prove Theorem~\ref{thm0.1}.
\begin{cor}
\label{cor5.10}
Let $F\colon \c C\to A$ be a connected Galois covering with group
$G$. Let $T=T_1\oplus\ldots\oplus T_n$ ($T_i\in ind(A)$) and
$T'=T'_1\oplus\ldots \oplus T'_n$ ($T'_i\in ind(A)$) be basic tilting $A$-modules lying in a
same connected component of $\overrightarrow{\c K}_A$. Assume that
$T,T'\in mod_1(A)$ and fix isomorphisms $\lambda_i\colon
F_{\lambda}\widehat{T}_i\xrightarrow{\sim} T_i$ and $\lambda_i\colon
F_{\lambda}\widehat{T'}_i\xrightarrow{\sim} T'_i$ with
$\widehat{T}_i,\widehat{T'}_i\in ind(\c C)$ for every $i$.
Then:
\begin{center}
$F_{\widehat{T}_i,\lambda_i}$ is connected $\Leftrightarrow$
$F_{\widehat{T'}_i,\lambda'_i}$ is connected.
\end{center}
In particular, if $T'=A$ or $T'=DA$, then
$F_{\widehat{T}_i,\lambda_i}$ is connected.
\end{cor}
\noindent{\textbf{Proof}}:
Recall that ``the Galois covering $\c E\to\c B$ is connected'' means $\c E$ is
connected and locally bounded. Thanks to Remark~\ref{rem2.5} and to
Proposition~\ref{prop2.11}, we know that $F_{\widehat{T}_i,\lambda_i}$
(resp. $F_{\widehat{T'}_i,\lambda'_i}$) is connected if and only if
$F.T$ (resp. $F.T'$) is a basic $\c C$-module. The corollary is
therefore a consequence of Proposition~\ref{prop5.9}.\hfill$\square$\\

We finish with the last objective of the section: give ``simple''
conditions on $T\in\overrightarrow{\c K}_A$ under which condition
($H_{A,T}$) is satisfied (see Lemma~\ref{lem2.3}).
\begin{prop}
\label{prop5.2}
Let $T,T'$ lie in a same connected component of $\overrightarrow{\c
K}_A$. Then:
\begin{center}
($H_{A,T}$) is satisfied $\Leftrightarrow$ ($H_{A,T'}$) is satisfied
\end{center}
In particular, if $T'=A$ or $T'=DA$, then ($H_{A,T}$) is satisfied.
\end{prop}
\noindent{\textbf{Proof:}} Let $\psi\colon A\xrightarrow{\sim} A$ be
an automorphism restricting to the identity map on objects. Let $x\in
Ob(A)$ and let $_?A_x\colon y\in Ob(A)\mapsto\ _yA_x$ be the
indecomposable projective $A$-module associated to $x$. Then, we have
an isomorphism of $A$-modules:
\begin{equation}
\begin{array}{rcl}
_?A_x & \longrightarrow & \psi.\,_?A_x\\
u\in\ _yA_x&\longmapsto & F(u)\in \left(\psi.\,_?A_x\right)(y)=\,_yA_x
\end{array}\notag
\end{equation}
So ($H_{A,A}$) is satisfied, and similarly ($H_{A,DA}$) is satisfied. Therefore, in order to prove the
proposition, it suffices to prove that ($H_{A,T}$) is satisfied if and
only if ($H_{A,T'}$) is satisfied, for any arrow $T\to T'$ in
$\overrightarrow{\c K}_A$. Assume that $T\to T'$ is an arrow in
$\overrightarrow{\c K}_A$ and that ($H_{A,T}$) is satisfied. So we
have the data:
\begin{enumerate}
\item[-] $T=X\oplus\overline{T}$ with $X\in ind(A)$,
\item[-] $T'=Y\oplus\overline{T}$ with $Y\in ind(A)$,
\item[-] a non split exact sequence $0\to
X\xrightarrow{u}M\to Y\to 0$ where $M\in add(\overline{T})$
and where $u\colon X\to M$ is the left
$add(\overline{T})$-approximation of $X$.
\end{enumerate}
Notice that in order to prove that ($H_{A,T'}$) is satisfied, we only
need to prove that $\psi.Y\simeq Y$. Since $\psi$ is an automorphism,
$\psi.\colon mod(A)\to mod(A)$ is an equivalence of abelian
categories. Therefore, the sequence $0\to
\psi.X\xrightarrow{\psi.u}\psi.M\to\psi.Y\to 0$ is non split exact and
verifies: $\psi.M\in add(\psi.\overline{T})$ and
$\psi.u\colon\psi.X\to\psi.M$ is the left
$add(\psi.\overline{T})$-approximation of $\psi.Y$. Moreover
$\psi.X\simeq X$,
$\psi.M\simeq M$ and $\psi.\overline{T}\simeq
\overline{T}$ because ($H_{A,T}$) is satisfied. So, $\psi.Y$ is
isomorphic to the cokernel of the left
$add(\overline{T})$-approximation of $X$. This implies that $Y\simeq
\psi.Y$. So ($H_{A,T'}$) is satisfied. The converse is dealt with
using dual arguments.\hfill$\square$\\

\section{Comparison of  $\protect\overrightarrow{\c K}_A$ and
  $\protect\overrightarrow{\c K}_{End_A(T)}$ for a tilting $A$-module $T$}\label{section4}

Let $T$ be a basic tilting $A$-module. Let
$B=End_A(T)$. 
In the preceding section, we have pointed out conditions of the form:
``$T$ lies in the connected component of $\overrightarrow{\c K}_A$
containing $A$''. Since our final objective (i.e. to compare the
Galois coverings of $A$ and of $B$) is symmetrical between $A$ and $B$,
we ought to find sufficient conditions for $T$ to lie in both
connected components of $\overrightarrow{\c K}_A$ and
$\overrightarrow{\c K}_B$ containing $A$ and $B$ respectively. Thus,
this section is devoted to compare $\overrightarrow{\c K}_A$ and
$\overrightarrow{\c K}_B$.
 For simplicity, if $X\in mod(A)$ (resp. $u\in Hom_A(X,Y)$) we
shall write  $X_T$ (resp. $u_T$) for the $B$-module (resp. the morphism
of $B$-modules) $Hom_A(X,T)$ (resp. $Hom_A(u,T)\colon Hom_A(Y,T)\to
Hom_A(X,T)$). Also, whenever $f$ is a morphism of modules, we shall
write $f_*$ (resp. $f^*$) for the mapping $g\mapsto fg$
(resp. $g\mapsto gf$).
We begin with a useful lemma.
\begin{lem}
\label{lem4.1}
Let $X\in mod(A)$ and let $T'\in \overrightarrow{\c K}_A$ be a predecessor
of $T$ (i.e. there is an oriented path in $\overrightarrow{\c K}_A$ starting at
$T'$ and ending at $T$). Then, the is an isomorphism, for any $Y\in add(T')$:
\begin{equation}
\begin{array}{crcl}
\theta_{X,Y}\colon & Hom_A(X,Y) & \longrightarrow & Hom_B(Y_T,X_T)\\
&u&\longmapsto & u_T
\end{array}\notag
\end{equation}
In particular: $Y\in
ind(A)\Leftrightarrow Y_T\in ind(B)$, for any $Y\in add(T')$.
\end{lem}
\noindent{\textbf{Proof:}} Remark that $\theta_{X,T'}$ is an
isomorphism if and only if $\theta_{X,Y}$ is an isomorphism for any
$Y\in add(T')$.
By assumption on
$T'$, there exists a path in $\overrightarrow{\c K}_A$ starting at $T'$
and ending at $T$. Let us prove by induction on the length $l$ of this
path that $\theta_{X,T'}$ is an isomorphism.

If $l=0$ then $T=T'$. So $\theta_{X,T'}=\theta_{X,T}$ is equal to:
\begin{equation}
\begin{array}{rcl}
Hom_A(X,T)=X_T&\longrightarrow & Hom_B(T_T,X_T)=Hom_B(B,X_T)\\
u& \longmapsto&(f\mapsto fu)
\end{array}\notag
\end{equation}
So $\theta_{X,T'}$ is an isomorphism (with inverse $\varphi\mapsto
\varphi(1_B)$). This proves the lemma when $l=0$.

Now, assume that $l>0$ and assume that $\theta_{X,T''}$ is an
isomorphism whenever $T''$ is the source of a path in
$\overrightarrow{\c K}_A$ ending at $T$ and with length equal to
$l-1$. We have a path $T'\to T''\to\ldots\to T$ of length $l$ in
$\overrightarrow{\c K}_A$. Therefore:
\begin{equation}
\text{$\theta_{X,Y}$ is an isomorphism for any $Y\in add(T'')$}\tag{$i$}
\end{equation}
Moreover, thanks to the arrow $T'\to T''$ in $\overrightarrow{\c K}_A$, we
have:
\begin{enumerate}
\item[$(ii)$] $T'=\overline{T}\bigoplus Y'$ with $Y'\in ind(A)$,
\item[$(iii)$] $T''=\overline{T}\bigoplus Y''$ with $Y''\in ind(A)$,
\item[$(iv)$] a non split exact sequence $0\to Y'\to M\to Y''\to 0$
with $M\in add(\overline{T})$.
\end{enumerate}
Thanks to $(i)$, $(ii)$ and $(iii)$ we only need to prove that
$\theta_{X,Y'}$ is an isomorphism. Remark that by assumption on $T'$
and $T''$ we have $T\in T^{\perp}\subseteq T''^{\perp}$. This implies
in particular that $Ext_A^1(Y'',T)=0$. Therefore, $(iv)$ yields an exact
sequence in $mod(A)$:
\begin{equation}
0\to Y''_T\to M_T\to Y'_T\to 0\notag
\end{equation}
This gives rise to the exact sequence:
\begin{equation}
0\to Hom_B(Y'_T,X_T)\to Hom_B(M_T,X_T)\to
Hom_B(Y''_T,X_T)\notag
\end{equation}
On the other hand, $(iv)$ yields the following exact sequence:
\begin{equation}
0\to Hom_A(X,Y')\to Hom_A(X,M)\to Hom_A(X,Y'')\notag
\end{equation}
Therefore, we have a commutative diagram:
\begin{equation}
\xymatrix{
0\ar@{->}[r] & Hom_A(X,Y')\ar@{->}[r] \ar@{->}[d]_{\theta_{X,Y'}} &
Hom_A(X,M) \ar@{->}[r]
\ar@{->}[d]_{\theta_{X,M}} & Hom_A(X,Y'')
\ar@{->}[d]_{\theta_{X,Y''}}\\
0\ar@{->}[r] & Hom_A(Y'_T,X_T)\ar@{->}[r] &
Hom_A(M_T,X_T) \ar@{->}[r]
& Hom_A(Y''_T,X_T)
}\notag
\end{equation}
where the rows are exact and where $\theta_{X,M}$ and
$\theta_{X,Y''}$ are isomorphisms. This shows that $\theta_{X,Y'}$ is
an isomorphism. So
$\theta_{X,T'}$ is an isomorphism and the induction is finished. This
proves the first assertion of the lemma.
The second assertion is due to the functoriality of
$\theta_{X,Y}$.\hfill$\square$\\

\begin{rem}
\label{rem4.2}
Assume that $A$ is hereditary. Then Lemma~\ref{lem4.1} still
  holds if one replaces the hypothesis "$T'$ is a predecessor of $T$"
  by "$T'\geqslant T$" (i.e. $T^{\perp}\subseteq T'^{\perp}$). The
  proof is then a classical application of left $add(T)$-approximations.
\end{rem}

The following proposition is the base of the link between
$\overrightarrow{\c K}_A$ and $\overrightarrow{\c K}_B$: it explains
how to associate suitable tilting $B$-modules with tilting
$A$-modules.
\begin{prop}
\label{prop4.3}
Let $X\to Y$ be an arrow in $\overrightarrow{\c K}_A$ where $X$ and $Y$
are predecessors of $T$. Then:
\begin{equation}
X_T\in \overrightarrow{\c K}_B\Leftrightarrow Y_T\in\overrightarrow{\c
K}_A\notag
\end{equation}
If the two conditions of the above equivalence are
satisfied, then there is an arrow $Y_T\to X_T$ in $\overrightarrow{\c K}_B$.
\end{prop}
\noindent{\textbf{Proof:}} Let us assume that $Y_T\in
\overrightarrow{\c K}_B$ and let us show that $X_T\in\overrightarrow{\c
K}_B$ and that there is an arrow $Y_T\to X_T$ in $\overrightarrow{\c
K}_B$ (the proof of the remaining implication is then obtained by
exchanging the roles of $X$ and $Y$). The arrow $X\to Y$ in
$\overrightarrow{\c K}_A$ gives the following data:
\begin{enumerate}
\item[.] $X=M\bigoplus \overline{X}$ with $M\in ind(A)$,
\item[.] $Y=N\bigoplus\overline{X}$ with $N\in ind(A)$,
\item[.] $\e\colon 0\to M\xrightarrow{i} X'\xrightarrow{p} N\to 0$ is
a non split exact sequence in $mod(A)$ with $X'\in add(\overline{X})$.
\end{enumerate}
The tilting $A$-module $Y$ is a predecessor of $T$.
Hence $T\in T^{\perp}\subseteq Y^{\perp}$ and therefore
$Ext_A^1(N,T)=0$. We infer that $Hom_A(\e, T)$ gives an exact sequence
in $mod(B)$:
\begin{equation}
0\to N_T\xrightarrow{p_T}X'_T\xrightarrow{i_T}M_T\to 0\tag{$\e_T$}
\end{equation}
Notice that we also have:
\begin{enumerate}
\item[.] $X_T=M_T\bigoplus \overline{X}_T$,
\item[.] $Y_T=N_T\bigoplus\overline{X}_T$,
\item[.] $X_T'\in add(\overline{X}_T)$.
\end{enumerate}
Therefore, in order to prove that $X_T\in \overrightarrow{\c K}_B$ and
that there is an arrow $Y_T\to X_T$ in $\overrightarrow{\c K}_B$, we
only need to prove the following facts:
\begin{enumerate}
\item[$1)$] $\e_T$ does not split,
\item[$2)$] $M_T\in ind(B)$ and $N_T\in ind(B)$,
\item[$3)$] $pd_B(X_T)<\infty$,
\item[$4)$] $X_T$ is selforthogonal,
\item[$5)$] $X_T$ is the direct sum of $n$
  indecomposable $A$-modules and $X_T$ is basic.
\end{enumerate}
\indent{$1)$} Let us prove that $\e_T$ does not split. If $\e_T$
splits, then $i_T$ is a retraction:
\begin{equation}
(\exists\lambda\in Hom_B(M_T,X'_T))\ \ Id_{M_T}=i_T\circ \lambda\notag
\end{equation}
Since $M$ is a direct summand of $X\in \overrightarrow{\c K}_A$ and
since $X$ is a predecessor of $T$, Lemma~\ref{lem4.1} implies that
$\lambda=\pi_T$ with $\pi\in Hom_A(X', M)$. Thus we have $(\pi\circ
i)_T=(Id_M)_T$. Using again Lemma~\ref{lem4.1} we deduce that
$\pi\circ i=Id_M$ which is impossible because $\e$ does not
split. So $\e_T$ does not split.

$2)$ Lemma~\ref{lem4.1} implies that $M_T,N_T\in ind(B)$.

$3)$ Since we assumed that $Y_T\in \overrightarrow{\c K}_B$, we have
$pd_B(\overline{X}_T)<\infty$, $pd_B(X'_T)<\infty$ and
$pd_B(N_T)<\infty$. Hence $\e_T$ gives $pd_B(M_T)<\infty$. So
$pd_B(X_T)<\infty$.

$4)$ Let us prove that $X_T$ is selforthogonal. In this purpose,
we will use the following lemma.
\begin{lem}
\label{lem4.4}
Let $L\in add(X)$. Then,
the following morphism induced by $p_T\colon N_T\to X'_T$ is
surjective:
\begin{equation}
\begin{array}{crcl}
(p_T)^*&Hom_B(X'_T,L_T)&\longrightarrow &
Hom_B(N_T,L_T)\\
&f&\longmapsto & f\circ p_T
\end{array}\notag
\notag
\end{equation}
\end{lem}
\noindent{\textbf{Proof:}}  Since
$L\in add(X)$ and since $X\in\overrightarrow{\c K}_A$, we have
$Ext_A^1(L,M)=0$. Hence, $Hom_A(L,\e)$ gives rise
to a surjective morphism induced by $p$:
\begin{equation}
\begin{array}{rrcl}
p_*\colon & Hom_A(L,X') & \twoheadrightarrow &
Hom_A(L,N)\\
&f&\longmapsto & p\circ f
\end{array}\notag
\end{equation}
Let us apply Lemma~\ref{lem4.1} to $X'\in
add(Y)$ and to $N\in add(Y)$. We get the following commutative diagram
where vertical arrows are isomorphisms:
\begin{equation}
\xymatrix{
Hom_A(L,X') \ar@{->}[r]^{p_*}
\ar@{->}[d]_{\theta_{L,X'}} & Hom_A(L,N)
\ar@{->}[d]_{\theta_{L,N}}\\
Hom_B(X'_T,L_T) \ar@{->}[r]^{(p_T)^*} &
Hom_B(N_T,L_T)
}\notag
\end{equation}
Since $p_*$ is surjective, we infer that so is
$(p_T)^*$.\hfill$\square$\\

Now we can prove that $X_T=\overline{X}_T\bigoplus M_T$ is
selforthogonal. Since $\overline{X}_T\in add(Y_T)$ and $Y_T\in \overrightarrow{\c
K}_B$, we get:
\begin{equation}
(\forall i\geqslant 1)\ \ Ext_B^i(\overline{X}_T,\overline{X}_T)=0\tag{$i$}
\end{equation}
For each $i\geqslant 1$, $Hom_B(\overline{X}_T,\e_T)$ gives the
following exact sequence:
\begin{equation}
Ext_B^i(\overline{X}_T,X'_T)\to Ext_B^i(\overline{X}_T,M_T)\to
Ext_B^{i+1}(\overline{X}_T,N_T)\notag
\end{equation}
Since $\overline{X}_T,X'_T,N_T\in add(Y_T)$ and $Y_T\in
\overrightarrow{\c K}_B$, we get:
\begin{equation}
(\forall i\geqslant 1)\ \ Ext_B^i(\overline{X}_T,M_T)=0\tag{$ii$}
\end{equation}
On the other hand, $Hom_B(\e_T,\overline{X}_T)$ gives the following
exact sequences:
\begin{align}
&.Hom_B(X'_T,\overline{X}_T)\xrightarrow{(p_T)^*}
Hom_B(N_T,\overline{X}_T)\to Ext_B^1(M_T,\overline{X}_T)\to
Ext_B^1(X'_T,\overline{X}_T)\notag\\
&.Ext_B^i(N_T,\overline{X}_T)\to Ext_B^{i+1}(M_T,\overline{X}_T)\to
Ext_B^{i+1}(X'_T,\overline{X}_T)\ \  \text{for $i\geqslant 1$}\notag
\end{align}
These exact sequences together with Lemma~\ref{lem4.4} and the
selforthogonality of $Y_T$ imply that:
\begin{equation}
(\forall i\geqslant 1)\ \ Ext_B^i(M_T,\overline{X}_T)=0\tag{$iii$}
\end{equation}
In order to get the selforthogonality of $X_T=M_T\bigoplus
\overline{X}_T$ it only remains to prove that $M_T$ is selforthogonal
(because of $(i),(ii)$ and $(iii)$). Notice that $Hom_B(N_T,\e_T)$
gives the following exact sequence for each $i\geqslant 1$:
\begin{equation}
Ext_B^i(N_T,X'_T)\to Ext_B^i(N_T,M_T)\to Ext_B^{i+1}(N_T,N_T)\notag
\end{equation}
Using $Y_T\in\overrightarrow{\c K}_B$ and $X'_T,N_T\in add(Y_T)$ we deduce that:
\begin{equation}
(\forall i\geqslant 1)\ \ Ext_B^i(N_T,M_T)=0\tag{$iv$}
\end{equation}
Finally $Hom_B(\e_T,M_T)$ gives the following exact sequences:
\begin{align}
&.Hom_B(X'_T,M_T)\xrightarrow{(p_T)^*}Hom_B(N_T,M_T)\to
Ext_B^1(M_T,M_T)\to Ext_B^1(X'_T,M_T)\notag\\
&.Ext_B^i(N_T,M_T)\to Ext_B^{i+1}(M_T,M_T)\to
Ext_B^{i+1}(X'_T,M_T)\  \ \text{for $i\geqslant 1$}\notag
\end{align}
These exact sequences together with Lemma~\ref{lem4.4}, $(ii)$ and
$(iv)$ imply that (recall that $X_T'\in add(\overline{X}_T)$):
\begin{equation}
(\forall i\geqslant 1)\ \ Ext_B^i(M_T,M_T)=0\tag{$v$}
\end{equation}
From $(i)$, $(ii)$, $(iii)$ and $(v)$ we deduce that $X_T=M_T\bigoplus
\overline{X}_T$ is selforthogonal.

$5)$ To finish, let us prove that $X_T$ is basic and that
$X_T$ is the direct sum of $n$ indecomposable modules. Notice that $\overline{X}_T$ is basic because it is a
direct summand of the basic tilting $B$-module $Y_T$.
On the other hand, $\e_T$ does not split, so $Ext_A^1(M_T,N_T)\neq
0$, hence $M_T\not\in add(Y_T)$ and therefore $M_T\not\in
add(\overline{X}_T)$. Since $M_T\in ind(B)$, we deduce that $X_T$ is basic.
Finally, $Y_T$ is by assumption the direct sum of $n$ indecomposable
modules, and $X_T$ and $Y_T$ differ by one indecomposable direct
summand so $X_T$ is also the direct sum of $n$ indecomposable modules.
\hfill$\square$\\

\begin{rem}
\label{rem4.7}
When $A$ is
hereditary, Proposition~\ref{prop4.3} has the following
generalisation: \textit{Let $X\in\overrightarrow{\c K}_A$ be such that
$X\geqslant T$, then $X_T\in\overrightarrow{\c K}_B$}. The proof of
this generalisation is obtained by replacing the use of the exact
sequence $\e$ by a coresolution of $X$ in $add(T)$. 
\end{rem}

Proposition~\ref{prop4.3} gives the following proposition
which will be used in the comparison of the Galois coverings of $A$
and $B$. We omit the proof which is immediate using Proposition~\ref{prop4.3}.
\begin{prop}
\label{prop4.6}
Let $X\in \overrightarrow{\c K}_A$ be such that there exists a path in
$\overrightarrow{\c K}_A$ starting at $X$ and ending at $T$. Then
$X_T\in\overrightarrow{\c K}_B$ and there exists in
$\overrightarrow{\c K}_B$ a path starting at $B$ and ending at
$X_T$.
\end{prop}

Proposition~\ref{prop4.3} also allows us to prove the main result of
this section. Recall that for a quiver $Q$, we write $Q^{op}$ for the
opposite quiver (obtained from $Q$ by reversing the arrows).
\begin{thm}
\label{thm4.5}
Let $\overrightarrow{\c K}_A(T)$ (resp. $\overrightarrow{\c K}_B(T)$)
be the convex hull of $\{A,T\}$
(resp. $\{B,T\}$) in $\overrightarrow{\c K}_A$
(resp. $\overrightarrow{\c K}_B$). Then we have an isomorphism of quivers:
\begin{equation}
\begin{array}{crcl}
\alpha\colon & \overrightarrow{\c K}_A(T)&\longmapsto&\overrightarrow{\c K}_B(T)^{op}\\
&X&\longmapsto & X_T=Hom_A(X,T)
\end{array}\notag
\end{equation}
 Under this correspondence, $A\in\overrightarrow{\c K}_A(T)$
 (resp. $T\in \overrightarrow{\c K}_A(T)$) is associated with $T\in\overrightarrow{\c K}_B(T)$
 (resp. $B\in\overrightarrow{\c K}_B(T)$).
\end{thm}
\noindent{\textbf{Proof:}} Thanks to Proposition~\ref{prop4.3}, the
mapping $\alpha$ is a well defined morphism of quivers. Thus, it only
remains to exhibit an inverse morphism. Notice that Proposition~\ref{prop4.6}
 implies that $\overrightarrow{\c K}_A(T)=\{A,T\}\
 \Leftrightarrow\overrightarrow{\c K}_B(T)=\{B,T\}\ \Leftrightarrow$
 there is no path in $\overrightarrow{\c K}_A(T)$ starting at $A$ and
 ending at $T$. Therefore, we may assume that there is a path starting
 at $A$ and ending at $T$. This assumption implies that any $X\in
 \overrightarrow{\c K}_A(T)$ is a predecessor of $T$.
From \cite[Thm 1.5]{miya} we
know that $T$ is a basic tilting $End_B(T)$-module and that we
have an isomorphism of $k$-algebras:
\begin{equation}
\begin{array}{rcl}
A&\longrightarrow&End_B(T)\\
a & \longmapsto & (t\mapsto at)
\end{array}\notag
\end{equation}
Henceforth, we shall consider $A$-modules as $End_B(T)$-modules and
vice-versa using the above isomorphism. In particular, we have an
identification of quivers:
\begin{equation}
\begin{array}{rcl}
\overrightarrow{\c K}_A(T) & \xrightarrow{\sim} & \overrightarrow{\c
  K}_{End_B(T)}(T)\\
X & \mapsto & X
\end{array}\notag
\end{equation}
Therefore, we also have a well defined morphism of quivers:
\begin{equation}
\begin{array}{crcl}
\alpha'\colon& \overrightarrow{\c K}_B(T)^{op}&\to &\overrightarrow{\c
  K}_A(T)\\
& X &\mapsto& X_T=Hom_B(X,T)
\end{array}\notag
\end{equation}
Let us prove that $\alpha'\alpha$
is an isomorphism. Let $X\in\overrightarrow{\c K}_A(T)$. Then $X$ is a
predecessor of $T$. Therefore, Lemma~\ref{lem4.1}
  implies that:
\begin{equation}
Hom_B(Hom_A(X,T),T)\simeq Hom_B(Hom_A(X,T),Hom_A(A,T))\simeq
Hom_A(A,X)\simeq X\notag
\end{equation}
This proves that $\alpha'\alpha$ is an isomorphism of
quivers. With the same arguments one also shows that $\alpha\alpha'$
is an isomorphism. So does $\alpha\colon \overrightarrow{\c
  K}_A(T)\to\overrightarrow{\c K}_B(T)^{op}$.\hfill$\square$\\

Notice that $\overrightarrow{\c
K}_A$ and $\overrightarrow{\c K}_B^{op}$ are not isomorphic in general. Indeed these quivers may have
different number of vertices as the following example shows.
\begin{ex}
\label{ex4.9}
Let $Q$ be the quiver:
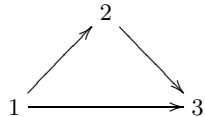
\begin{equation}
\xymatrix{
&2\ar@{->}[rd]&\\
1\ar@{->}[rr] \ar@{->}[ru] && 3
}\notag
\end{equation}
and let $A=kQ/I$ where $I$ is the ideal generated by the oriented path
of length $2$ in $Q$. Notice that $A$ is of finite representation
type. Let $T=P_1\oplus P_2\oplus \tau_A^{-1}P_3$ be
the APR-tilting $A$-module associated to the sink $3$. Hence:
\begin{equation}
T=
\begin{matrix}
1\\
2\ 3
\end{matrix} \oplus
\begin{matrix}
2\\3
\end{matrix} \oplus
\begin{matrix}
\ \ 1\ \ 2\\
2\ \ 3
\end{matrix}\notag
\end{equation}
and the Hasse diagram $\overrightarrow{\c K}_A$ of basic tilting $A$-modules is equal to:
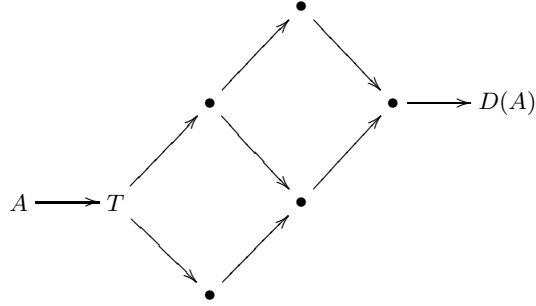
\begin{equation}
\xymatrix{
&&&\bullet\ar@{->}[rd] &&\\
&& \bullet\ar@{->}[rd] \ar@{->}[ru]&& \bullet\ar@{->}[r] &D(A)\\
A \ar@{->}[r] & T \ar@{->}[ru] \ar@{->}[rd] && \bullet\ar@{->}[ru]&&\\
&&\bullet\ar@{->}[ru] &&&
}\notag
\end{equation}
On the other hand, $B=End_A(T)$ is isomorphic to $kQ'/I'$ where $Q'$
is equal to:
\begin{equation}
\xymatrix{
a \ar@{->}[r] & b \ar@/^/@{->}[r] & c\ar@/^/@{->}[l]
}\notag
\end{equation}
and $I'$ is the ideal generated by the path $c\to b\to c$. As a
$B$-module, $T$ is equal to
\begin{equation}
T=
\begin{matrix}
a\\b\\c\\b
\end{matrix}\oplus
\begin{matrix}
c\\b
\end{matrix}\oplus
\begin{matrix}
a\ c\\
b
\end{matrix}\notag
\end{equation}
and $\overrightarrow{\c K}_B$ is equal to:
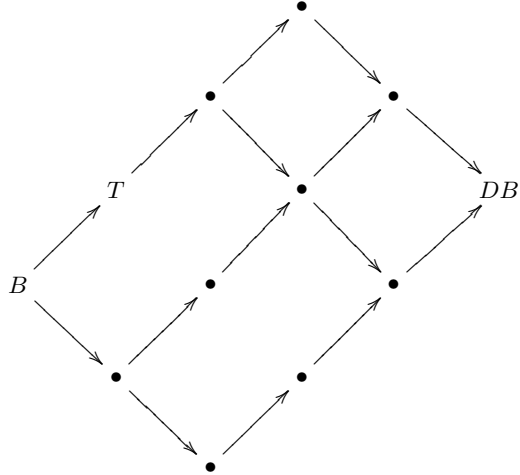
\begin{equation}
\xymatrix{
&&&\bullet\ar@{->}[rd]&&\\
&&\bullet\ar@{->}[rd]\ar@{->}[ru]&&\bullet\ar@{->}[rd]&\\
&T\ar@{->}[ru]&&\bullet\ar@{->}[ru]\ar@{->}[rd]&&DB\\
B\ar@{->}[ru]\ar@{->}[rd]&&\bullet\ar@{->}[ru]&&\bullet\ar@{->}[ru]&\\
&\bullet\ar@{->}[rd]\ar@{->}[ru]&&\bullet\ar@{->}[ru]&&\\
&&\bullet\ar@{->}[ru]&&&
}\notag
\end{equation}
In particular, $\overrightarrow{\c K}_A$ and $\overrightarrow{\c
K}_B$ do not have the same number of vertices. Notice that, in this example, the isomorphism of
Theorem~\ref{thm4.5} is equal to:
\begin{equation}
\begin{array}{rcl}
\overrightarrow{\c K}_A(T)=\left(A\to T\right) &\longrightarrow & 
\overrightarrow{\c K}_B(T)^{op}=\left(T\to B\right)\\
A &\longmapsto & T=Hom_A(A,T)\\
T & \longmapsto & B=Hom_A(T,T)
\end{array}\notag
\end{equation}
\end{ex}

\begin{rem}
\label{rem4.8}
Assume that $A$ is hereditary, then Theorem~\ref{thm4.5} has the
following generalisation, thanks to Remark~\ref{rem4.7}: \textit{Let
  $Q_A$ (resp. $Q_B$) be the full subquiver of $\overrightarrow{\c
    K}_A$ (resp. $\overrightarrow{\c K}_B$) made of the tilting
  modules $X\geqslant T$. Then $X\mapsto X_T$ induces an isomorphism
  of quivers $Q_A\xrightarrow{\sim} Q_B^{op}$}.
\end{rem}

\section{Comparison of the Galois coverings of $A$ and $End_A(T)$ for
  $T$ basic tilting $A$-module}\label{section5}

This section is devoted to the proof of Theorem~\ref{thm0.1}, of
Corollary~\ref{cor0.3} and of Corollary~\ref{cor0.2}. Let
$T\in \overrightarrow{\c K}_A$ and let $B=End_A(T)$. 
As in the introduction, we shall say that $A$ and $B$ have the same
connected Galois coverings with group $G$ if and only if there exists
a bijection $Gal_A(G)\xrightarrow{\sim} Gal_B(G)$. Here $Gal_A(G)$
denotes the set of equivalence classes of connected Galois coverings
with group $G$ of $A$.
In order to
compare the equivalence classes of connected Galois coverings of $A$
and those of $B$, we introduce the following assertion which depends
on $A$, on $T$ and on a fixed group $G$:
\begin{center}
$\c P(A,T,G)=$''\textit{($H_{A,T}$) is satisfied and for any connected Galois covering $F\colon \c
    C\to A$ with group $G$, the $A$-module $T$ is of the first kind
    w.r.t. $F$ and $F.T$ is a basic $\c C$-module}''
\end{center}
Recall from Definition~\ref{def2.4} that the condition ($H_{A,T}$)
ensures the existence of an equivalence class $[F]_T$ of Galois
coverings of $B$ depending only on the equivalence class $[F]$ of
$F$. Recall also from Remark~\ref{rem2.5} and from
Proposition~\ref{prop2.11} that the condition ``$F.T$ is a basic $\c
C$-module'' implies that $[F]_T$ is an equivalence class of connected
Galois coverings of $B$.
Finally, recall that $\c P(A,A,G)$ and $\c P(A,DA,G)$ are true for any
$G$ (see Proposition~\ref{prop5.1}, Proposition~\ref{prop5.9} and
Proposition~\ref{prop5.2}). The above definition of 
$\c P(A,T,G)$ is relevant because of the following proposition.
\begin{prop}
\label{prop6.1}
Let $G$ be a group. Assume that $\c P(A,T,G)$ and $\c P(B,T,G)$ are
true. Then $A$ and $B$ have the same connected Galois coverings with
group $G$.
\end{prop}
\noindent{\textbf{Proof}}: Since $\c P(A,T,G)$ is true,  we
have a well defined mapping:
\begin{equation}
\begin{array}{rrcl}
\varphi_A\colon & Gal_A(G) & \longrightarrow & Gal_B(G)\\
&[F] & \longmapsto & [F]_T
\end{array}\tag{$i$}
\end{equation}
Similarly, $\c P(B,T,G)$ is true so we have a well defined mapping:
\begin{equation}
\begin{array}{rrcl}
\varphi_B\colon & Gal_B(G) & \longrightarrow & Gal_G(End_B(T))\\
&[F] & \longmapsto & [F]_T
\end{array}\tag{$ii$}
\end{equation}
Thanks to  Proposition~\ref{prop2.7} we know that $\rho_A^{-1}\circ
 \left(\varphi_B\varphi_A([F])\right)=[F]$ for any $[F]\in Gal_A(G)$. 
Therefore, $\varphi_A$ is one-to-one and $\varphi_B$ is onto.
Notice that thanks to the isomorphism $\rho_A\colon
 A\xrightarrow{\sim} End_B(T)$, the assertion $\c P(End_B(T),T,G)$ is
 true, so that the above arguments imply that $\varphi_B$ is one-to-one and
 that $\varphi_{End_B(T)}$ is onto.
As a consequence, $\varphi_B$ is bijective, so the mapping $[F]\mapsto
 \rho_A^{-1}\circ [F]_T$ induces a bijection $Gal_B(G)\xrightarrow{\sim}Gal_A(G)$.
\hfill$\square$\\

Thanks to Proposition~\ref{prop6.1} we are reduced to find sufficient
conditions for $\c P(A,T,G)$ and $\c P(B,T,G)$ to be simultaneously
true.
The following
proposition is a direct consequence of Proposition~\ref{prop5.6},
 of Corollary~\ref{cor5.10}, of Proposition~\ref{prop5.2} and of the fact that
$\c P(A,A,G)$ and $\c P(A,DA,G)$ are true.
\begin{prop}
\label{prop6.3}
Let $G$ be a group. Let $T'\in \overrightarrow{\c K}_A$ lying in the
connected component of $\overrightarrow{\c K}_A$ containing $T$. Then:
\begin{center}
$\c P(A,T,G)$ is true $\Leftrightarrow$ $\c P(A,T',G)$ is true
\end{center}
In particular, if $T'=A$ or if $T'=DA$ then $\c P(A,T,G)$ is true.
\end{prop}

Thanks to Proposition~\ref{prop6.3}, we are reduced look for conditions for
$T$ to lie in both connected components of $\overrightarrow{\c K}_A$
and $\overrightarrow{\c K}_B$ containing $A$ and $B$
respectively. Such a condition is given by the following proposition.
\begin{prop}
\label{prop6.4}
Let $G$ be a group and assume that there exists a path in
$\overrightarrow{\c K}_A$ starting at $A$ and ending at $T$. Then $T$
lies in the connected component of $\overrightarrow{\c K}_A$
(resp. $\overrightarrow{\c K}_B$) containing $A$
(resp. $B$). Consequently, $\c
P(A,T,G)$ and $\c P(B,T,G)$ are true.
\end{prop}
\noindent{\textbf{Proof}}: Theorem~\ref{thm4.5} implies that there
exists a path in $\overrightarrow{\c K}_B$ starting at $Hom_A(T,T)=B$ and ending
at $Hom_A(A,T)=T$. Using Proposition~\ref{prop6.3} we get the desired
conclusion.\hfill$\square$\\

Now we can prove Theorem~\ref{thm0.1}:\\
\noindent{\textbf{Proof of Theorem~\ref{thm0.1}}}:
$1)$ Since $T$ and $T'$ lie in a same connected component of
$\overrightarrow{\c K}_A$, there exists a sequence $T^{(1)}=T,T^{(2)},\ldots,T^{(r)}=T'$ of
basic tilting $A$-modules such that for any $i\in\{1,\ldots,r-1\}$,
there exists a path in $\overrightarrow{\c K}_A$ with $T^{(i)}$ and
$T^{(i+1)}$ as end-points. For short, let us write $B_i$ for
$End_A(T^{(i)})$. Let $i\in\{1,\ldots,r-1\}$ and let us
assume, for example, that there exists a path in $\overrightarrow{\c
  K}_A$ starting at $T^{(i)}$ and ending at $T^{(i+1)}$. Using Lemma~\ref{lem4.1}
and Proposition~\ref{prop4.6} we infer that:
\begin{enumerate}
\item[$(i)$] $End_A(T^{(i)})$ and $End_{B_{i+1}}(Hom_A(T^{(i)},T^{(i+1)}))$ are isomorphic as
  $k$-algebras (and therefore as $k$-categories),
\item[$(ii)$] there exists a path in $\overrightarrow{\c K}_{B_{i+1}}$
  starting at $B$ and ending at $Hom_A(T^{(i)},T^{(i+1)})$.
\end{enumerate}
This implies (thanks to Proposition~\ref{prop6.4} and to
Proposition~\ref{prop6.1}) that $End_A(T^{(i)})$
and $End_A(T^{(i+1)})$ have the same connected Galois coverings with group
$G$. Since this fact is true for any $i$, we deduce that $End_A(T)$
and $End_A(T')$ have the same connected Galois coverings with group $G$.

$2)$ is a consequence of $1)$, of the fact that $End_A(A)\simeq
End_A(DA)\simeq A^{op}$
and of the fact that $A$ and $A^{op}$ have the same Galois coverings
($F\colon\c C\to A$ is a Galois covering if and only if
$F^{op}\colon\c C^{op}\to A$ is a Galois covering and $\c C^{op}$ is
connected and locally bounded if and only if $\c
C^{op}$ is).\hfill$\square$\\

Using Theorem~\ref{thm0.1} we can prove Corollary~\ref{cor0.3} and
Corollary~\ref{cor0.2}.

\noindent{\textbf{Proof of Corollary~\ref{cor0.3}}:}
Since $A$ is of finite representation type, Theorem~\ref{thm0.1}
implies that $A$ and $B$ have the same connected Galois
coverings. For the same reason, $A$
(resp. $B$) admits a connected Galois covering with group $G$ if and
only if $G$ is a factor group of the fundamental group $\pi_1(A)$
(resp. $\pi_1(B)$) of
$A$ (resp. of $B$). Consequently, $\pi_1(A)$ and $\pi_1(B)$ are
isomorphic.\hfill$\square$\\ 

\noindent{\textbf{Proof of Corollary~\ref{cor0.2}}}: $1)$ and $2)$ are consequences of
Theorem~\ref{thm0.1} and of the fact that $A$ is simply connected if and
only if it has no proper connected Galois covering (see
\cite[Cor. $4$]{lemeur}).

$3)$ is a consequence of $2)$.
\hfill$\square$\\

Corollary~\ref{cor0.3} naturally leads to the following question: let
$G$ be a group such that $A$ and $B$ have the same Galois coverings
with group $G$, is it true that $A$ admits an admissible presentation
with fundamental group isomorphic to $G$ if and only if the same holds
for $B$? The answer is no in general as the following example shows :
\begin{ex}
\label{ex6.5}
Let $Q$ be the following quiver:
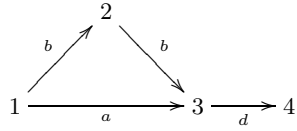
\begin{equation}
\xymatrix{
&2\ar@{->}[rd]^b &&\\
1\ar@{->}[rr]_a \ar@{->}[ru]^b & & 3\ar@{->}[r]_d & 4
}\notag
\end{equation}
and let $A=kQ/I$ where $I=<da>$.
Let $T=P_1\oplus P_2\oplus P_3\oplus \tau_{A}^{-1}(P_4)=P_1\oplus
P_2\oplus P_3\oplus S_3$ be the APR-tilting module associated with the
sink $3$ (here $S_i$ is the simple $A$-module associated to the vertex
$i$ and $P_i$ is the indecomposable projective $A$-module with top
$S_i$). Then $B=End_A(T)$ is the path algebra of the following quiver:
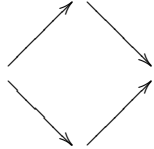
\begin{equation}
\xymatrix{
&\ar@{->}[rd]&\\
\ar@{->}[ru] \ar@{->}[rd] &&\\
&\ar@{->}[ru]
}\notag
\end{equation}
Since $T$ is an APR-tilting $A$-module, there is an arrow $A\to T$ in
$\overrightarrow{\c K}_A$. Then, Theorem~\ref{thm0.1} implies that for
any group $G$, the $k$-algebras $A$ and $B$ have the same connected
Galois covering with group $G$. On the other hand, any admissible
presentation of $B$ has fundamental group isomorphic to $\mathbb{Z}$
whereas $A$ admits an admissible presentation with fundamental group
$0$ and another one with fundamental group isomorphic to
$\mathbb{Z}$ (see for example \cite[1.4]{assem_delapena}).
\end{ex}

In the preceding example, the reader may remark that the
fundamental group of any admissible presentation of $A$ is a factor
group of $\mathbb{Z}$ and that the same holds for $B$. Let us say that
$A$ admits an optimum fundamental group ($G$) if and only if there
exists an admissible presentation of $A$ with fundamental group $G$
and if the fundamental group of any other admissible presentation is a
factor group of $G$. For example, $A$ admits an optimum fundamental
group in the following cases: $A$ is of finite representation type
(see \cite{gabriel}), $A$ is constricted (see \cite[Thm 3.5]{bardzell_marcos}),
$A$ is monomial, $A$ is triangular and has no double bypass (see
\cite[Thm. 1]{lemeur}). Then we have the following corollary whose proof is a
direct consequence of Theorem~\ref{thm0.1}:
\begin{cor}
Assume that $T$ lies in the connected component of $\overrightarrow{\c
K}_A$ containing $A$. Then $A$ admits $G$ as optimum fundamental group
if and only if $B$ admits $G$ as optimum fundamental group.
\end{cor}

\section{On the simple connectedness of a tilted algebra}\label{sections}

The aim of this section is to prove Proposition~\ref{prop0.4}. We
shall use the construction of the preceding
section. Recall from Remark~\ref{rem2.5}, Lemma~\ref{lem2.9} and
Proposition~\ref{prop2.11} that given a tilting
$A$-module $T$ and a connected Galois covering $F\colon \c C\to A$, we
need two properties in order to define a connected Galois
covering of $End_A(T)$: 
\begin{enumerate}
\item $T$ is of the first w.r.t. $F$,
\item  $F.T$ is a
basic $\c C$-module.
\end{enumerate}
For the first property, we use the following
lemma which is due to W.~Crawley-Boevey (see
\cite{geiss_leclerc_schroer}). I acknowledge B.~Keller and C.~Geiss
who gave me the details about this lemma.
\begin{lem}
\label{lem7.1}(W.~Crawley-Boevey, \cite{geiss_leclerc_schroer})
Let $G=\mathbb{Z}/p\mathbb{Z}$ where $p$ is the characteristic of
$k$. ($G=\mathbb{Z}$ if $car(k)=0$). Let $\c E,\c B$ be locally
bounded $k$-categories
and let $F\colon \c E\to \c B$ be a
Galois covering with group $G$. Let
$M$ be a rigid indecomposable
$\c B$-module. Then $M$ is of the first kind w.r.t. $F$.
\end{lem}

\begin{rem}
\label{rem7.2}
Lemma~\ref{lem7.1} also holds for any group
in the class of groups containing $\mathbb{Z}/p\mathbb{Z}$ and closed
under direct product, under quotient and under isormophism.
\end{rem}

\begin{lem}
\label{lem7.3}
Let $G$ be a finite group, let $F\colon A'\to A$ be a connected Galois
covering with group $G$ and let $T$ be a basic tilting $A$-module of
the first kind w.r.t. $F$. Then $F.T$ is a basic $A'$-module.
\end{lem}
\noindent{\textbf{Proof:}} First, let us prove that $F.T$ is a basic
$A'$-module. From Lemma~\ref{lem1.1} applied to $T$, we have
$pd_{A'}(F.T)<\infty$. Besides, $F_{\lambda}F.T\simeq\oplus_{g\in G}T$
because $T$ is of the first kind w.r.t. $F$, so that:
\begin{equation}
(\forall i\geqslant 1)\ \ Ext_{A'}^i(F.T,F.T)\simeq
Ext_A^i(F_{\lambda}F.T,T)\simeq Ext_A^i(\bigoplus_{g\in G}T,T)\simeq
\bigoplus_{g\in G}Ext_A^i(T,T)=0\notag
\end{equation}
Thus, $F.T$ is an exceptionnal $A'$-module. Moreover, if we apply $F.$
to any finite coresolution of $A$ in $add(T)$, we deduce a finite
coresolution of $A'\simeq F.A$ in $add(F.T)$. Hence $F.T$ is a tilting $A'$-module.
Now we can prove that $F.T$ is basic. Since $T$ is basic, it is the
direct sum of $n$ indecomposable
$A$-modules. Moreover, $T$ is of the first kind w.r.t. $F$, so $F.T$
is the direct sum of $n.|G|=rk(K_0(A'))$ indecomposable
$A'$-modules. Since $F.T$ is a tilting $A'$ module, this implies that
$F.T$ is basic.\hfill$\square$\\

Now we can prove Proposition~\ref{prop0.4}.\\
\noindent{\textbf{Proof of Proposition~\ref{prop0.4}:}} 
Recall (see for example \cite{happel2}) that $HH^1(kQ)=0$ if and only
if $Q$ is a tree.
Now assume
that $Q$ is not a tree. Then, $A$ admits a connected Galois covering
$F\colon \c C\to A$ with group
$\mathbb{Z}$. Let $G$ be any finite cyclic group if $car(k)=0$ and
let $G=\mathbb{Z}/p\mathbb{Z}$ if $car(k)=p$ is non zero. Let $N$ be the
kernel of the natural surjection $\mathbb{Z}\twoheadrightarrow G$. Then there
exists a commutative diagram:
\begin{equation}
\xymatrix{
\c C\ar@{->}[dd]_F \ar@{->}[rd]^p&\\
&\c C/N\ar@{->}[ld]^{\overline{F}}\\
A &
}\notag
\end{equation}
where $p\colon \c C\to \c C/N$ is the natural projection and
$\overline{F}\colon \c C/N\to A$ is a connected Galois covering with
group $G$. Then:
\begin{enumerate}
\item[-] if $car(k)=0$, then Lemma~\ref{lem7.1} implies that $T$ is of
  the first kind w.r.t. $F$ so
  that it is also of the first kind w.r.t. $\overline{F}$.
\item[-] if $car(k)=p$ is non zero then Lemma~\ref{lem7.1} implies
  that $T$ is of the first kind w.r.t. $\overline{F}$.
\end{enumerate}
 Thanks to Lemma~\ref{lem7.3} we deduce that $\overline{F}\colon\c
 C/N\to A$ is
 connected Galois covering of $A$ with group $G$ and verifying: $T$ is
 of the first kind w.r.t. $F$ and $F.T$ is a basic $\c
 C/N$-module. Using Remark~\ref{rem2.5} and
Proposition~\ref{prop2.11}, we deduce that $End_A(T)$ admits a
 connected Galois covering with group $G\neq 1$. Hence $End_A(T)$ is
 not simply connected.\hfill$\square$\\

\section*{Final remark}
The Hasse diagram $\overrightarrow{\c K}_A$ of basic tilting
$A$-modules describes the combinatoric relations between tilting
modules. When $A$ is hereditary (i.e. $A=kQ$ with $Q$ a finite quiver
with no oriented cycle) these combinatorics are also described by
the cluster category $\c C_Q$ of the quiver $Q$ (see \cite{bmrrt}). In
particular, the indecomposable tilting objects in $\c C_Q$ are
displayed as the vertices of an unoriented graph. Since this graph is
always connected (see \cite[3.5]{bmrrt}) it is natural to
ask if it is possible to remove all conditions concerning connected components
in Theorem~\ref{thm0.1} and Corollary~\ref{cor0.2} (in the hereditary
case). These
 developpements will be detailed in a forecoming text.
\bibliographystyle{plain}
\bibliography{biblio}

\begin{thebibliography}{10}

\bibitem{assem_skowronski}
I.~Assem and Skowro\'nski A.
\newblock On some classes of simply connected algebras.
\newblock {\em Proceedings of the {L}ondon {M}athematical {S}ociety},
  56(3):417--450, 1988.

\bibitem{assem_skowronski2}
I.~Assem and Skowro\'nski A.
\newblock Tilting simply connected algebras.
\newblock {\em Communications in {A}lgebra}, 22(12):4611--4619, 1994.

\bibitem{assem_coelho_trepode}
I.~Assem, F.~U. Coelho, and S.~Trepode.
\newblock Simply connected tame quasi-tilted algebras.
\newblock {\em Journal of {P}ure and {A}pplied {A}lgebra}, 172(2--3):139--160,
  2002.

\bibitem{assem_delapena}
I.~Assem and J.~A. de~La~Pe\~na.
\newblock The fundamental groups of a triangular algebra.
\newblock {\em Communications in {A}lgebra}, 24(1):187--208, 1996.

\bibitem{assem_marcos_delapena}
I.~Assem, E.~N. Marcos, and J.~A. de~La~Pe\~na.
\newblock The simple connectedness of a tame tilted algebra.
\newblock {\em Journal of {A}lgebra}, 237(2):647--656, 2001.

\bibitem{apr}
M.~Auslander, M.~I. Platzeck, and I.~Reiten.
\newblock Coxeter functors without diagrams.
\newblock {\em Transactions of the {A}merican {M}athematical {S}ociety},
  250:1--46, 1979.

\bibitem{auslander_reiten}
M.~Auslander and I.~Reiten.
\newblock Applications of contravariantly finite subcategories.
\newblock {\em Advances in {M}athematics}, 86(1):111--152, 1991.

\bibitem{bardzell_marcos}
M.~J. Bardzell and E.~N. Marcos.
\newblock ${H}^1$ and representation of finite dimensional algebras.
\newblock {\em Lecture {N}otes in {P}ure and {A}pplied {M}athematics},
  224:31--38, 2002.

\bibitem{bongartz_gabriel}
K.~Bongartz and P.~Gabriel.
\newblock Covering spaces in representation theory.
\newblock {\em Inventiones {M}athematicae}, 65:331--378, 1982.

\bibitem{brenner_butler}
S.~Brenner and M.~C.~R. Butler.
\newblock Generalisation of the {B}ernstein-{G}elfand-{P}onomarev reflection
  functor.
\newblock {\em Lecture {N}otes in {M}athematics}, 832:103--169, 1980.

\bibitem{bmrrt}
A.B. Buan, R.~Marsh, M.~Reineke, I.~Reiten, and G.~Todorov.
\newblock Tilting theory and cluster combinatorics.
\newblock math.RT/0402054, 2004.

\bibitem{gabriel}
P.~Gabriel.
\newblock The universal cover of a representation finite algebra.
\newblock {\em Lecture {N}otes in {M}athematics}, 903:65--105, 1981.
\newblock in: Representation of algebras.

\bibitem{geiss_leclerc_schroer}
C.~Geiss, B.~Leclerc, and J.~Schroër.
\newblock Auslander algebras and initial seeds for cluster algebras.
\newblock math.RT/0506405, 2005.

\bibitem{happel}
D.~Happel.
\newblock {\em Triangulated categories in the representation theory of finite
  dimensional algebras}, volume 119 of {\em London {M}athematical {S}ociety
  {L}ecture {N}otes {S}eries}.
\newblock Cambridge {U}niversity {P}ress, {C}ambridge, 1988.

\bibitem{happel2}
D.~Happel.
\newblock Hochschild cohomology of finite dimensional algebras.
\newblock {\em Séminaire d'{A}lgèbre {P}aul {D}ubreuil, {M}arie-{P}aule
  {M}alliavin, {L}ecture {N}otes in {M}athematics}, 1404:108--126, 1989.

\bibitem{happel_ringel}
D.~Happel and C.~M. Ringel.
\newblock Tilted algebras.
\newblock {\em Transactions of the {A}merican {M}athematical {S}ociety},
  274(2):399--443, 1982.

\bibitem{happel_unger}
D.~Happel and L.~Unger.
\newblock On a partial order of tilting modules.
\newblock {\em Algebras and {R}epresentation {T}heory}, 8(2):147--156, 2005.

\bibitem{happel_unger2}
D.~Happel and L.~Unger.
\newblock On the quiver of tilting modules.
\newblock {\em Journal of {A}lgebra}, 284(2):857--868, 2005.

\bibitem{lemeur}
P.~Le~Meur.
\newblock The universal cover of an algebra without double bypasses.
\newblock math.RT/0507513, 2005.

\bibitem{martinezvilla_delapena}
R.~Mart\'inez-Villa and J.~A. de~la Pe\~na.
\newblock The universal cover of a quiver with relations.
\newblock {\em Journal of {P}ure and {A}pplied {A}lgebra}, 30:277--292, 1983.

\bibitem{miya}
Y.~Miyashita.
\newblock Tilting modules of finite projective dimension.
\newblock {\em Mathematische {Z}eitschrift}, 193(1):113--146, 1986.

\bibitem{riedtmann}
C.~Riedtmann.
\newblock Algebren, darstellungsköcher, Überlagerungen und zurück.
\newblock {\em Commentarii {M}athematici {H}elvetici}, 55(2):199--224, 1980.

\bibitem{riedtmann_schofield}
C.~Riedtmann and A.~Schofield.
\newblock On a simplicial complex associated with tilting modules.
\newblock {\em Commentarii {M}athematici {H}elvetici}, 66(1):70--78, 1991.

\bibitem{skowronski2}
A.~Skowro{\'n}ski.
\newblock Simply connected algebras and {H}ochschild cohomologies.
\newblock {\em Canadian {M}athematical {S}ociety {C}onference {P}roceedings},
  14:431--447, 1993.

\end{thebibliography}
\end{document}